\documentclass[11pt]{article}
\usepackage{graphicx}
\usepackage{color}
\usepackage{amsmath,amssymb,amsthm,amsfonts}
\usepackage{amssymb}
\usepackage{cite}

\usepackage{subfig}
\usepackage{amsmath}
\allowdisplaybreaks[4]
\newtheorem{thm}{Theorem}[section]
\newtheorem{cor}[thm]{Corollary}
\newtheorem{lem}[thm]{Lemma}
\newtheorem{prop}[thm]{Proposition}

\theoremstyle{definition}
\newtheorem{defn}{Definition}[section]
\theoremstyle{remark}

\numberwithin{equation}{section}

\DeclareMathSymbol{\C}{\mathalpha}{AMSb}{"43}

\newcommand{\alp}{\alpha}

\newcommand{\R}{{\mathbb{R}}}

\newcommand{\bsub}{\begin{subequations}}
\newcommand{\esub}{\end{subequations}$\!$}

\begin{document}
\title{Asymptotic Behavior of the Principal Eigenvalue Problems with Large Divergence-Free Drifts}

\author{ Yujin Guo\thanks{School of Mathematics and Statistics, and  Key Laboratory of Nonlinear Analysis $\&$ Applications (Ministry of Education), Central China Normal University, P.O. Box 71010, Wuhan 430079, P. R. China. Email: \texttt{yguo@ccnu.edu.cn}.}
,\, 	Yuan Lou\thanks{School of Mathematical Sciences, CMA-Shanghai and MOE-LSC, Shanghai Jiao Tong University, Shanghai 200240,
 P. R.	China. Email: \texttt{yuanlou@sjtu.edu.cn}.},\,
and\, Hongfei Zhang\thanks{School of Mathematics and Statistics, and  Hubei Key Laboratory of
Mathematical Sciences, Central China Normal University, P.O. Box 71010, Wuhan 430079,
P. R. China. Email: \texttt{hfzhang@mails.ccnu.edu.cn}.}
}
\renewcommand{\thefootnote}{}
\date{\today }
\smallbreak\maketitle
\begin{abstract}
In this paper, we consider the following principal  eigenvalue problem with a large divergence-free drift:
\begin{equation}\label{0.1}
-\varepsilon\Delta \phi-2\alpha\nabla m(x)\cdot\nabla \phi+V(x)\phi=\lambda_\alpha \phi\ \,\ \text{in}\, \ H_0^1(\Omega),
\end{equation}
where the domain $\Omega\subset \R^N (N\ge 1)$ is bounded with smooth boundary $\partial\Omega$, the constants $\varepsilon>0$ and $\alpha>0$ are the diffusion and drift coefficients, respectively, and $m(x)\in C^{2}(\bar{\Omega})$, $V (x)\in C^{\gamma}(\bar{\Omega})~(0<\gamma<1)$ are given functions.
For a class of divergence-free drifts where $m$ is a harmonic function
in $\Omega$ and has no first integral in $H_{0}^{1}(\Omega)$, we prove the convergence of the principal eigenpair $(\lambda_\alpha, \phi)$ for (\ref{0.1}) as $\alpha\rightarrow+\infty$, which addresses a special case of the open question proposed in [H.  Berestycki, F. Hamel and N.
Nadirashvili, CMP, 2005]. Moreover, we further investigate the refined limiting profiles of the principal eigenpair  $(\lambda_\alpha, \phi)$ for (\ref{0.1}) as $\alpha\rightarrow+\infty$, which display the visible effects of the large divergence-free drifts on the principal eigenpair $(\lambda_\alpha, \phi)$.
\end{abstract}

\vskip 0.05truein


	
\vskip 0.1truein

\section{Introduction}

In this paper, we are concerned with the following principal  eigenvalue problem with a   drift:
\begin{equation}\label{1.1}
\left\{
\begin{split}
-\varepsilon\Delta \phi-2\alpha\nabla m(x)\cdot\nabla \phi+V(x)\phi&=\lambda_\alpha \phi\ \ \text{in}\ \ \Omega,\\
\phi&=0\  \,\ \quad\hbox{on}\, \ \partial\Omega,
\end{split}
\right.
\end{equation}
where the origin is an interior point of the bounded domain $\Omega\subset\mathbb{R}^{N}~(N\geq1)$ with smooth boundary $\partial\Omega$, and
$\varepsilon>0$, $\alpha>0$ represent the
diffusion and drift coefficients, respectively. Here the drift term $\nabla m$ satisfying $m(x)\in C^{2}(\bar{\Omega})$ is {\em divergence-free, $i.e.$, $\Delta m\equiv0$ in $\Omega$}, and $V(x)$ satisfies the general condition $0\leq V (x)\in C^{\gamma}(\bar{\Omega})~(0<\gamma<1)$. The principal eigenvalue problem (\ref{1.1}) has strong motivations in  the nonlinear propagation phenomena, the competition of species, and the stability analysis of equilibria, see \cite{BHN,CC2,DEF,DF,Friedman} and the references therein. Particularly, Friedman and his collaborators \cite{DEF,DF,Friedman} studied in 1970s the estimates of the principal eigenvalue and eigenfunction for (\ref{1.1}).
We also remark that the principal eigenvalue problem (\ref{1.1}) with other boundary conditions was previously studied in \cite{CL,CL1,PZZ,PZ}, and the parabolic version of (1.1) was explored in \cite{LLPZ,LLPZ1} and the references therein.

It is known from \cite{Guo,PZ,PZZ} that the problem (\ref{1.1}) admits a unique principal ($i. e.$, first) eigenvalue $\lambda_{\alpha}$, which corresponds to a unique (up to the multiplication) principal ($i. e.$, first) eigenfunction  $\phi_{\alpha}>0$. Throughout the present whole paper, for convenience, the above unique pair $(\lambda_\alpha, \phi_\alpha)$ is referred to as the unique principal eigenpair of (\ref{1.1}).
When the drift term $\nabla m$ is not divergence-free, the limiting profiles of the unique principal eigenpair for (\ref{1.1}) were recently studied  in our previous work \cite{Guo} as $\alpha\rightarrow+\infty$. However, when the drift term $\nabla m$ is divergence-free, the problem (\ref{1.1}) is more relevant to fluid dynamics and mathematical biology, which was addressed widely over the past few years, see  \cite{BHN,LL2} and references therein. Especially, starting from the pioneering work of Berestycki,  Hamel and Nadirashvili in \cite{BHN}, the monotonicity and other analytical properties of the principal eigenvalue $\lambda_{\alpha}$ for (\ref{1.1}) were investigated in \cite{GGP1,GGP2,IKNR,MV,SSSZ} and the references therein.

In order to study the principal eigenvalue problem (\ref{1.1}) with a  divergence-free drift, the existing works  (e.g. \cite{BHN}) often introduce the following definition of first integrals:

\begin{defn}\label{def1}
A nonzero function $\omega\in H^{1}(\Omega)$ is said to be a first integral of the vector field $\nu$, if $\nu\cdot\nabla \omega=0$ holds a.e. in $\Omega$. We denote $\mathcal{I}_{0}$ the set of all first integrals for the vector field $\nu$, which belong to $  H_{0}^{1}(\Omega)$.
\end{defn}

\noindent More precisely, when the drift term $\nabla m$ is divergence-free and has first integrals in $H_{0}^{1}(\Omega)$, the convergence of the principal eigenfunction $\phi_{\alpha}$ for (\ref{1.1}) was proved  in \cite{BHN}. However, when the drift term $\nabla m$ is divergence-free and has no first integral in $H_{0}^{1}(\Omega)$, the convergence of the  principal eigenfunctions for (\ref{1.1}) was proposed as an open question in \cite[Section 5]{BHN}. Of course, one may expect from \cite{BHN,Guo} that the answer of this interesting question depends usually on the specific forms of $m(x)$ and $V(x)$.

Stimulated by above facts, the purpose of the present paper is to investigate the asymptotic behavior of the unique principal eigenpair for (\ref{1.1}) as $\alpha\rightarrow+\infty$, provided that the drift term $\nabla m$ is divergence-free and has no first integral in $H_{0}^{1}(\Omega)$.
Towards this purpose, we consider the divergence-free case where $m(x)$ satisfies for $x=(x_1,\cdots,x_N)\in\bar{\Omega}$,
\begin{equation}\label{1.4}
 m(x)=a_{0}+\sum_{i=1}^{N}a_{i}x_{i}^{2},\ \text{and}\ \, a_{0}\in\mathbb{R},\ a_{i}\in\mathbb{R}\setminus \{0\}\, \
 \text{satisfies}\, \ \sum_{i=1}^{N}a_{i}=0.
\end{equation}
Under the assumption (\ref{1.4}), one can check that the drift term $\nabla m$ is divergence-free, $i.e.$, $\Delta m\equiv0$ in $\Omega$, and has no first integral in $H_{0}^{1}(\Omega)$, see Lemma \ref{lemB} below for details.
We also note from \cite{Guo,PZZ} that the principal eigenvalue $\lambda_{\alpha}$ of $(\ref{1.1})$ can be  characterized equivalently by
\begin{equation}\label{1.02}
\begin{split}
\lambda_{\alpha}=&\inf_{\{\phi\in H_{0}^{1}(\Omega),~\phi\not\equiv 0\}}
\frac{\int_{\Omega}e^{\frac{2\alpha}{\varepsilon}m(x)}(\varepsilon|\nabla \phi|^{2}+V(x)\phi^{2})dx}{\int_{\Omega}e^{\frac{2\alpha}{\varepsilon}m(x)}\phi^{2}dx}\\
=&\inf_{\{u\in H_{0}^{1}(\Omega),~\int_{\Omega}u^{2}dx=1\}}\Big\{\int_{\Omega}\varepsilon|\nabla u|^2dx+\int_{\Omega}\Big[\frac{\alpha^{2}}{\varepsilon}|\nabla m|^{2}
+\alpha \Delta m
+V(x)\Big]u^2dx\Big\},
\end{split}
\end{equation}
where $\varepsilon>0$, $\alpha>0$ and $u=e^{\frac{\alpha }{\varepsilon}m}\phi$.
Using the variational theory, if $u_{\alpha}>0$ is a positive minimizer of (\ref{1.02}), then $(\lambda_{\alpha}, u_{\alpha})$ satisfies the following elliptic equation
\begin{equation}\label{1.31}
-\varepsilon\Delta u_{\alpha}+\Big[\frac{\alpha^{2}}{\varepsilon}|\nabla m|^{2}+\alpha \Delta m+V(x)\Big]u_{\alpha}
=\lambda_{\alpha}u_{\alpha}
\ \ \text{in}\ \, H_{0}^{1}(\Omega).
\end{equation}
Therefore, suppose $(\lambda_{\alpha}, \phi_{\alpha})$ is the unique principal eigenpair of (\ref{1.1}),  throughout the present whole paper, it always means that $(\lambda_{\alpha}, u_{\alpha})$ satisfies both (\ref{1.02}) and (\ref{1.31}), where $u_\alpha>0$ is normalized in the sense that $\|u_\alpha\|^2_{L^2(\Omega)}=1$.

Under the assumption (\ref{1.4}), the main result of the present paper is concerned with the following convergence.

\begin{thm}\label{thm1.2}
Suppose  $m(x)$ satisfies (\ref{1.4}), and assume $0\leq V (x)\in C^{\gamma}(\bar{\Omega})~(0<\gamma<1)$. Then for any fixed $\varepsilon>0$, the unique principal eigenpair $(\lambda_{\alpha}, u_{\alpha})$ of (\ref{1.1}) satisfies
\begin{equation}\label{1.5}
\lim_{\alpha\rightarrow+\infty}\alpha^{-1}\lambda_{\alpha}=2\sum_{j=1}^{N}|a_{j}|>0,
\end{equation}
and
\begin{equation}\label{1.6}
\begin{split}
\hat{w}_\alpha(x)&:=\alpha^{-\frac{N}{4}}u_{\alpha}\big(\alpha^{-\frac{1}{2}}x+d_\alpha\big)\rightarrow Q(x):=\Big(\prod_{i=1}^{N}\frac{2|a_{i}|}{\varepsilon\pi}\Big)^{\frac{1}{4}}
e^{-\frac{1}{\varepsilon}\sum_{j=1}^{N}|a_{j}|x_{j}^{2}}>0\\
& \ \ \text{strongly in}\ \ H^{1}(\mathbb{R}^{N})\cap L^{\infty}(\mathbb{R}^{N})
\ \ \text{as}\ \ \alpha\rightarrow+\infty,
\end{split}
\end{equation}
where $a_j\in \R\setminus\{0\}$ is as in (\ref{1.4}) for $j=1,\cdots,N$, $u_{\alpha}(x)\equiv0$ in $\R^{N}\backslash\Omega$, and $d_\alpha\in\Omega$ is a global maximum point of $u_\alpha$ and satisfies
\begin{equation}\label{2.0021m}
\lim_{\alpha\to +\infty}\alpha^{\frac{1}{2}}|d_\alpha|=0.
\end{equation}
\end{thm}

Suppose $m(x)$ satisfies the assumption (\ref{1.4}), $i.e.,$ the case where  $\nabla m$ is divergence-free and has no first integral in $H_{0}^{1}(\Omega)$, Theorem \ref{thm1.2} gives
the explicit convergence of the unique principal eigenpair for (\ref{1.1}) as $\alpha\to +\infty$. This thus addresses a special case of the open question proposed in \cite[Section 5]{BHN} for an important class of divergence-free drifts. Moreover, one can note from (\ref{1.5}) and (\ref{1.6}) that the leading terms of the unique principal eigenpair $(\lambda_{\alpha},u_{\alpha})$ for (\ref{1.1})  as $\alp\rightarrow+\infty$ explicitly depend only on the drift term $\nabla m(x)$, instead of the trapping potential $V(x)$. Of course, it also implies from (\ref{1.6}) and (\ref{2.0021m}) that the mass of the principal eigenfunction $u_\alpha$ for (\ref{1.1}) concentrates near the unique saddle point ($i.e.,$ the origin) of $m(x)$ as $\alpha\rightarrow+\infty$. Furthermore, we shall analyze in Corollary \ref{cor2.3} that the global maximum point $d_\alpha\in\Omega$ of the principal eigenfunction $u_\alpha$ for (\ref{1.1}) must be unique as $\alpha\rightarrow+\infty$, provided that $|a_i|=|a_j|>0$  in (\ref{1.4}) holds for all  $i\neq j$. One may further wonder whether the global maximum point $d_\alpha\in\Omega$ of the principal eigenfunction $u_\alpha$ for (\ref{1.1}) is still unique as $\alpha\rightarrow+\infty$, particularly in the case where $3\le N\in \mathbb{N}^+$ is odd and $|a_i|\not =|a_j|$ in (\ref{1.4}) holds for some $i\neq j$.


The proof of Theorem \ref{thm1.2} makes full use of the following limiting variational problem
\begin{equation}\label{AA:A1}
\mu:=\inf_{\{u\in S,~\int_{\mathbb{R}^{N}}u^{2}dx=1\}}
\int_{\mathbb{R}^{N}}\Big[
\varepsilon|\nabla u|^{2}+\frac{4}{\varepsilon}\Big(\sum_{i=1}^{N}a_{i}^{2}x_{i}^{2}\Big)u^{2}\Big]dx,
\end{equation}
where $\varepsilon>0$ is the same as in Theorem \ref{thm1.2}, $a_{i}\in \R\setminus\{0\}$  is as in (\ref{1.4}) for $i=1,\cdots,N$, and the space $S$ is defined by
\begin{equation}\label{A2}
S:=\Big\{u\in H^{1}(\R^{N}):\,\ \int_{\mathbb{R}^{N}}\Big(\sum_{i=1}^{N}a_{i}^{2}x_{i}^{2}\Big)u^{2}dx<\infty\Big\},
\end{equation}
together with the norm
$$\|u\|_{S}=\Big\{\int_{\R^{N}}\Big[|\nabla u|^{2}+\Big(1+\sum_{i=1}^{N}a_{i}^{2}x_{i}^{2}\Big)u^{2}\Big]dx\Big\}^{\frac{1}{2}}.$$
We shall prove in Section 2 that the minimum value $\mu$ of (\ref{AA:A1}) satisfies
\[
\mu=2\sum_{i=1}^{N}|a_{i}|>0,
\]
and (\ref{AA:A1}) admits a unique positive minimizer $Q(x)>0$ satisfying (\ref{1.6}), where   $\int_{\R^N}Q^2dx=1$,  and $a_i\in \R\setminus\{0\}$ is as in (\ref{1.4}) for $i=1,\cdots,N$. By studying the variational problem (\ref{AA:A1}), we shall derive finally the limiting energy (\ref{1.5}) of (\ref{1.1}) as $\alpha\rightarrow+\infty$. Employing the blow up analysis, together with the elliptic regularity theory,  in Section 2 we further obtain the explicit convergence  (\ref{1.6}) of the unique principal eigenfunction $u_\alpha$ for (\ref{1.1}) as $\alpha\to +\infty$.

Under the assumptions of Theorem \ref{thm1.2}, one can obtain from Theorem \ref{thm1.2} that for any fixed $\varepsilon>0$, the unique principal eigenfunction $u_{\alpha}$ of (\ref{1.1}) satisfies
\begin{equation}\label{M:1.6}
\begin{split}
&w_\alpha(x):=\alpha^{-\frac{N}{4}}u_{\alpha}\big(\alpha^{-\frac{1}{2}}x\big)\rightarrow Q(x)>0\\
  & \text{strongly  in}\, \ H^{1}(\mathbb{R}^{N})\cap L^{\infty}(\mathbb{R}^{N})
\ \ \text{as}\ \ \alpha\rightarrow+\infty,
\end{split}
\end{equation}
where  $Q(x)>0$ is given by (\ref{1.6}), and $u_{\alpha}(x)\equiv0$ in $\R^{N}\backslash\Omega$. Following the convergence (\ref{M:1.6}), we are further concerned with the refined limiting profiles of the unique principal eigenpair $(\lambda_{\alpha}, u_{\alpha})$ for (\ref{1.1}) as $\alpha\rightarrow+\infty$. We first denote $\varphi_{i}(x)\in C^{2}(\R^N)\cap L^{\infty}(\R^N)~(i=1,2)$ to be the unique solution of
\begin{equation}\label{1.7}
\left\{
\begin{split}
&-\varepsilon\Delta \varphi_{i}+\Big(\frac{4}{\varepsilon}\sum_{j=1}^{N}a_j^2x_j^2-2\sum_{j=1}^{N}|a_{j}| \Big)\varphi_{i}=F_{i}(x)\ \,\ \text{in}\ \,  \mathbb{R}^{N},\\
&\int_{\mathbb{R}^{N}}\varphi_{i}Qdx=0,\ \ i=1,2,
\end{split}\right.
\end{equation}
where  $a_{j}\in \R\setminus\{0\}$ is as in (\ref{1.4}) for $j=1,\cdots,N$, $Q(x)>0$ is given by (\ref{1.6}), and $F_{i}(x)$ satisfies
\begin{equation}\label{1.7m}
F_{i}(x)=
\left\{\begin{split}
\Big[\int_{\R^N}\big(x\cdot\nabla V(0)\big)Q^2dx-x\cdot\nabla V(0)\Big]Q(x),\quad  &\text{if}\ \ i=1;\qquad   \\
\sum_{|\tau|=2}\frac{D^{\tau}V(0)}{\tau!}\Big(\int_{\R^N}x^{\tau}Q^2dx-x^{\tau}\Big)Q(x), \qquad\ \
&\text{if}\ \ i=2.
\end{split}\right.
\end{equation}
Here $\tau=(\tau_{1},\cdots,\tau_{N})$ is a multiple index with nonnegative integers $\tau_{1},\cdots,\tau_{N}$, and
$$|\tau|:=\tau_{1}+\cdots+\tau_{N},\ \ \tau!:=\tau_{1}!\cdots\tau_{N}!,\ \  x^{\tau}:=x_{1}^{\tau_{1}}\cdots x_{N}^{\tau_{N}},$$
provided that $ x=(x_{1},\cdots,x_{N})\in \R^N$.
Since $\int_{\R^N}F_i(x)Q(x)dx=0$ holds for $i=1,2$ in view of $\int_{\R^N}Q^2dx=1$, similar to \cite[Lemma 3.4]{Guo}, one can obtain the existence and uniqueness of $\varphi_{i}(x)$ for $i=1,2$. Suppose $V(x)$  satisfies $V(x)\in C^{2}(B_{r_0}(0))$ in some small ball $B_{r_0}(0)\subset\Omega$, then we have the following refined limiting profiles of (\ref{1.1}) as $\alpha\rightarrow+\infty$.

\begin{thm}\label{thm1.3}
Suppose $m(x)$ satisfies (\ref{1.4}), and assume $0\leq V (x)\in C^{\gamma}(\bar{\Omega}) \cap C^{2}(B_{r_0}(0))$, where $0<\gamma<1$ and $B_{r_0}(0)\subseteq\Omega$ is a small ball for some $r_0>0$. Then for any fixed $\varepsilon>0$, the unique principal eigenpair $(\lambda_{\alpha}, u_{\alpha})$ of (\ref{1.1}) satisfies
\begin{equation}\label{1.8}
\begin{split}
\lambda_{\alpha}=&2\alpha\sum_{i=1}^{N}|a_{i}|+V(0)+\alpha^{-\frac{1}{2}}\int_{\R^N}[x\cdot \nabla V(0)]Q^2dx\\
&+\alp^{-1}\sum_{|\tau|=2}\frac{D^{\tau}V(0)}{\tau!}\int_{\R^N}x^{\tau}Q^2dx
+o(\alpha^{-1})\ \ \text{as}\ \ \alpha\rightarrow+\infty,
\end{split}
\end{equation}
and
\begin{equation}\label{1.9}
\begin{split}
w_{\alpha}(x):=&\alpha^{-\frac{N}{4}}u_{\alpha}
(\alpha^{-\frac{1}{2}}x)\\
=&Q(x)+\alpha^{-\frac{3}{2}}\varphi_1(x)+\alpha^{-2}\varphi_2(x)
+o(\alpha^{-2}) \ \ \text{in}\ \ \mathbb{R}^{N}
\ \, \text{as}\ \ \alpha\rightarrow+\infty,
\end{split}
\end{equation}
where $a_{i}\in \R\setminus \{0\}$ is as in (\ref{1.4}) for $i=1, \cdots, N$, $Q(x)>0$ is given by (\ref{1.6}),  $u_\alpha(x)\equiv0$ in $\mathbb{R}^{N}\backslash\Omega$, and $\varphi_1(x)$, $\varphi_2(x)$ are given uniquely by (\ref{1.7}) and (\ref{1.7m}).
\end{thm}

There are several comments on Theorem \ref{thm1.3} in order. Firstly, as far as we know, Theorem \ref{thm1.3} presents the first result on the refined limiting profiles of the unique principal eigenpair $(\lambda_{\alpha}, u_{\alpha})$ for (\ref{1.1}) as $\alp\rightarrow+\infty$, for the case where the drift term $\nabla m$ of (\ref{1.1}) is divergence-free. Secondly, one can note from Theorem \ref{thm1.3} that the form of $m(x)$ affects generally all order asymptotic  expansions of $(\lambda_{\alpha}, u_{\alpha})$ as $\alp\rightarrow+\infty$, while the potential $V(x)$ however affects $(\lambda_{\alpha}, u_{\alpha})$ starting from its  second-order asymptotic expansion as $\alp\rightarrow+\infty$. Thirdly, applying the convergence (\ref{1.5}), the refined limiting profile of $\lambda_{\alpha}$ is analyzed as $\alp\rightarrow+\infty$ by establishing the important Pohozaev type identity (\ref{m}). Applying  the convergence (\ref{M:1.6}) and the refined estimate of $\lambda_{\alpha}$ as  $\alp\rightarrow+\infty$, the refined limiting profile of $u_{\alpha}$ is further derived as $\alp\rightarrow+\infty$ by employing the $L^{2}$--constraint condition of $u_{\alpha}$.

One can note that the refined limiting profiles of Theorem \ref{thm1.3} require a higher regularity of the potential $V(x)$ near the origin. We next illustrate that the refined limiting profiles of $(\lambda_{\alpha}, u_{\alpha})$ can be still derived  as $\alp\rightarrow+\infty$, provided that the potential $V(x)$ has a lower regularity near the origin. Towards this purpose, we now introduce (cf. \cite{Grossi}) the following homogeneous functions.

\begin{defn}\label{def1.1}
A function $h(x): \R^N\mapsto\R$ is called homogeneous of degree
$k\in \R^{+}$ about the origin, if
\begin{equation}\label{v1}
h(tx)=t^{k}h(x)\ \ \mbox{for any}\ \, t\in\R^{+}\ \, \hbox{and}\ \, x\in\R^{N}.
\end{equation}
\end{defn}

\noindent
It yields from Definition \ref{def1.1} that if the function $h(x)\in C(\R^N)$ is homogeneous of degree $k>0$ about
the origin, then it follows that $|h(x)|\leq C|x|^{k}$ holds in $\R^{N}$, where $C=\max\limits_{x\in \partial B_{1}(0)}h(x)>0$.

We next denote $\varphi_{i}(x)\in C^{2}(\R^N)\cap L^{\infty}(\R^N)~(i=3,4)$ to be the unique solution of
\begin{equation}\label{1.10}
\left\{
\begin{split}
&-\varepsilon\Delta \varphi_{i}+\Big(\frac{4}{\varepsilon}\sum_{j=1}^{N}a_j^2x_j^2-2\sum_{j=1}^{N}|a_{j}| \Big)\varphi_{i}=F_{i}(x)\ \,\ \text{in}\ \,  \mathbb{R}^{N},\\
&\int_{\mathbb{R}^{N}}\varphi_{i}Qdx=0,\ \ i=3,4,
\end{split}\right.
\end{equation}
where $a_{j}\in \R\setminus \{0\}$ is as in (\ref{1.4}) for $j=1, \cdots, N$, $Q(x)>0$ is given by (\ref{1.6}), and $F_{i}(x)$ satisfies
\begin{equation}\label{1.11}
F_{i}(x)=
\left\{\begin{split}
\Big[\int_{\R^N}h_0(x)Q^2dx-h_0(x)\Big]Q(x),\quad\quad\qquad\qquad\ \ &\text{if}\ \ i=3;\\
\Big(\int_{\R^N}h_0(x)\varphi_3Qdx\Big)Q+\Big(\int_{\R^N}h_0(x)Q^{2}dx-h_0(x)\Big)\varphi_3,\ \,\
&\text{if}\ \ i=4.
\end{split}\right.
\end{equation}
Here $h_0(x)$ satisfies that
\begin{equation}\label{1.12}
\begin{split}
 0\leq h_0(x)&\in C^{\gamma}(\R^{N})~(0<\gamma<1)\ \, \text{is homogeneous of degree}\, \ k_0>0 \\
&\text{about the origin, and satisfies}\ \lim_{|x|\rightarrow+\infty}h_0(x)=+\infty.
\end{split}
\end{equation}
Similar to \cite[Lemma 3.4]{Guo}, one can also have the existence and uniqueness of $\varphi_i(x)$ for (\ref{1.10}) with $i=3, 4$. Suppose the regularity of the potential $V(x)$   near the origin is lower than that of Theorem \ref{thm1.3}, then the following theorem still gives the refined limiting profiles of $(\lambda_{\alpha}, u_{\alpha})$ as $\alp\rightarrow+\infty$.

\begin{thm}\label{cor1.3}
Suppose   $m(x)$ satisfies (\ref{1.4}), and assume   $0\leq V (x)\in C^{\gamma}(\bar{\Omega})~(0<\gamma<1)$ satisfies $V(x)=h_0(x)$ in some small ball $B_{\hat{r}_0}(0)\subseteq\Omega$ for some $\hat{r}_0>0$, where $h_0(x)$ satisfies (\ref{1.12}) for some $k_0>0$. Then for any fixed $\varepsilon>0$, the unique principal eigenpair $(\lambda_{\alpha}, u_{\alpha})$ of (\ref{1.1}) satisfies
\begin{equation}\label{1.13}
\begin{split}
\lambda_{\alpha}=&2\alpha\sum_{i=1}^{N}|a_{i}| +\alpha^{-\frac{k_0}{2}}\int_{\R^{N}}h_0(x)Q^2dx \\
&+\alpha^{-k_0-1}\int_{\R^{N}}h_0(x)\varphi_3Qdx+o(\alpha^{-k_0-1})\ \ \text{as}\ \ \alpha\rightarrow+\infty,
\end{split}
\end{equation}
and
\begin{equation}\label{1.14}
\begin{split}
w_{\alpha}(x):=&\alpha^{-\frac{N}{4}}u_{\alpha}
(\alpha^{-\frac{1}{2}}x)\\
=&Q(x)+\alpha^{-\frac{k_0+2}{2}}\varphi_3(x)+\alpha^{-k_0-2}\Big[\varphi_4(x)-\frac{1}{2}
\Big(\int_{\R^N}\varphi_3^2(x)dx\Big)Q(x)\Big]\\
&+o(\alpha^{-k_0-2}) \ \ \text{in}\ \ \mathbb{R}^{N}
\ \ \text{as}\ \ \alpha\rightarrow+\infty,
\end{split}
\end{equation}
where $a_{i}\in \R\setminus \{0\}$ is as in (\ref{1.4}) for $i=1, \cdots, N$, $Q(x)>0$ is given by (\ref{1.6}),  $u_\alpha(x)\equiv0$ in $\mathbb{R}^{N}\backslash\Omega$, and $\varphi_3(x)$, $\varphi_4(x)$ are given uniquely by (\ref{1.10}) and (\ref{1.11}).
\end{thm}

We emphasize that Theorem \ref{cor1.3} holds for any given degree $k_0>0$. Moreover, if one only needs the first two terms in the expansions of $(\lambda_{\alpha}, u_{\alpha})$ for  (\ref{1.1}) as $\alpha\rightarrow+\infty$, then the assumptions on $V(x)$ of Theorem \ref{cor1.3} can be relaxed, see Lemma \ref{lem4.1} below for more details. Furthermore, the proof of Theorem \ref{cor1.3} is similar to that of Theorem \ref{thm1.3}, for which one needs to make full use of the $L^{2}$--constraint condition for $u_{\alpha}$, and as well establish  an important Pohozaev type identity. We finally comment that our arguments can yield the similar results of Theorems 1.1--1.3, no matter whether the minimum points of $m(x)$ and $V(x)$ coincide or not.

This paper is organized as follows. In Section 2, we first analyze the limiting variational problem (\ref{AA:A1}), after which we shall prove Theorem \ref{thm1.2} in Subsection 2.1 on the convergence of the unique principal eigenpair $(\lambda_{\alpha}, u_{\alpha})$ for (\ref{1.1}) as $\alpha\rightarrow+\infty$. In Section 3, we first prove Lemmas \ref{lem3.1} and \ref{lem3.2} on the refined estimates of $(\lambda_{\alpha}, u_{\alpha})$ as $\alpha\rightarrow+\infty$, based on which we shall complete in Subsection 3.1 the proof of Theorem \ref{thm1.3}. The proof of Theorem \ref{cor1.3} is finally given in Section 4.

\section{Convergence of $(\lambda_{\alpha}, u_{\alpha})$ as $\alpha\rightarrow+\infty$}
The main purpose of this section is to prove Theorem \ref{thm1.2} on the convergence of the unique principal eigenpair $(\lambda_{\alpha}, u_{\alpha})$ for (\ref{1.1}) as $\alpha\rightarrow+\infty$, which thus addresses a special case of the open question proposed in \cite{BHN}.


Under the assumption (\ref{1.4}), which specially assumes that $\sum_{i=1}^{N}a_{i}=0$, one can check that  the drift field $\nabla m$ is divergence-free. The following lemma further shows that the drift field $\nabla m$ satisfying (\ref{1.4}) has no first integral in $H_0^1(\Omega)$.

\begin{lem}\label{lemB}
Suppose $m(x)$ satisfies (\ref{1.4}), then the drift field $\nabla m$ does not have any first integral in $H_{0}^{1}(\Omega)$.
\end{lem}

\noindent\textbf{Proof.}
On the contrary, suppose the drift field $\nabla m$ has a first integral $\omega(x)\in H_{0}^{1}(\Omega)$. It then follows from Definition \ref{def1} that $\omega(x)$ satisfies
\begin{equation}\label{B1}
\omega\not\equiv 0\ \ \text{and}\ \ \nabla m\cdot\nabla \omega=0\ \ a. e.\ \, \text{in}\ \, \Omega.
\end{equation}
Without loss of generality, one can further assume from (\ref{B1}) that $\omega(x)$ satisfies
\begin{equation}\label{B01}
\omega(x)\geq0\ \ a. e.\ \, \text{in}\ \, \Omega.
\end{equation}
Indeed, denote
$$\omega(x):=\omega^{+}(x)-\omega^{-}(x)\in H_{0}^{1}(\Omega),$$
where $\omega^{+}(x):=\max\{\omega(x), 0\}\in H_{0}^{1}(\Omega)$ and $\omega^{-}(x):=\max\{-\omega(x), 0\}\in H_{0}^{1}(\Omega)$. Since it follows from (\ref{B1}) that
$$\nabla m\cdot\nabla \omega=\nabla m\cdot\nabla \omega^{+}-\nabla m\cdot\nabla \omega^{-}=0\ \ a. e.\ \, \text{in}\ \, \Omega,$$
one can have
\begin{equation}\label{B2}
\nabla m\cdot\nabla \omega^{+}=0 \ \ \text{and}\ \ \nabla m\cdot\nabla \omega^{-}=0\ \ a. e.\ \, \text{in}\ \, \Omega.
\end{equation}
We then conclude from (\ref{B2}) that $|\omega|$ satisfies
\begin{equation*}
\begin{split}
\nabla m\cdot\nabla |\omega|&=\nabla m\cdot\nabla (\omega^{+}+\omega^{-})\\
&=\nabla m\cdot\nabla\omega^{+}+\nabla m\cdot\nabla\omega^{-}
=0\ \ a. e.\ \, \text{in}\ \, \Omega,
\end{split}
\end{equation*}
which further implies that $|\omega|\in H_{0}^{1}(\Omega)$ is also a first integral of the drift field $\nabla m$. This proves (\ref{B01}).

Since $m(x)$ satisfies (\ref{1.4}), we derive from (\ref{B1}) that
\begin{equation}\label{B3}
\nabla m(x)\cdot\nabla\omega(x)
=2\sum_{i=1}^{N}a_{i}x_{i}\omega_{x_{i}}=0\ \ a. e.\ \, \text{in}\ \, \Omega,
\end{equation}
where $x=(x_1,\cdots,x_N)\in\Omega$.
Multiplying the equation (\ref{B3}) by $x_{1}^{2}$ and integrating over $\Omega$, we then deduce from (\ref{1.4}) that
\begin{align}
0&=\int_{\Omega}\Big(2\sum_{i=1}^{N}a_{i}x_{i}\omega_{x_{i}}
\Big)x_{1}^{2}dx\nonumber\\
&=2\Big[a_{1}\int_{\Omega}x_{1}^{3}\omega_{x_{1}}dx
+\int_{\Omega}\Big(\sum_{j=2}^{N}a_{j}x_{j}\omega_{x_{j}}
\Big)x_{1}^{2}dx\Big]\nonumber\\
&=-2\Big[3a_{1}\int_{\Omega}x_{1}^{2}\omega dx+
\Big(\sum_{j=2}^{N}a_{j}\Big)\int_{\Omega}x_{1}^{2}\omega dx\Big]
=-4a_{1}\int_{\Omega}x_{1}^{2}\omega dx,\nonumber
\end{align}
which further implies from (\ref{B01}) that $\omega(x)\equiv 0$ $a. e.$ in $\Omega$. This is however a contradiction with (\ref{B1}), and the proof of Lemma \ref{lemB} is therefore complete.
\qed

In order to discuss the convergence of the unique principal eigenpair $(\lambda_{\alpha}, u_{\alpha})$ for (\ref{1.1}) as $\alpha\rightarrow+\infty$, we now introduce the following limiting variational problem
\begin{equation}\label{A1}
\mu:=\inf_{\{u\in S,~\int_{\mathbb{R}^{N}}u^{2}dx=1\}}
\int_{\mathbb{R}^{N}}\Big[
\varepsilon|\nabla u|^{2}+\frac{4}{\varepsilon}\Big(\sum_{i=1}^{N}a_{i}^{2}x_{i}^{2}\Big)u^{2}\Big]dx,
\end{equation}
where the fixed constant $\varepsilon>0$ is as in (\ref{1.1}),  $a_{i}\in \R\setminus\{0\}$ is as in (\ref{1.4}) for $i=1, \cdots, N$, and the space $S$ is defined by
\begin{equation}\label{A2}
S:=\Big\{u\in H^{1}(\R^{N}): \ \int_{\mathbb{R}^{N}}\Big(\sum_{i=1}^{N}a_{i}^{2}x_{i}^{2}\Big)u^{2}dx<\infty\Big\},
\end{equation}
together with the norm
$$\|u\|_{S}=\Big\{\int_{\R^{N}}\Big[|\nabla u|^{2}+\Big(1+\sum_{i=1}^{N}a_{i}^{2}x_{i}^{2}\Big)u^{2}\Big]dx\Big\}^{\frac{1}{2}}.$$
The following lemma describes  the minimum energy $\mu>0$ of the problem (\ref{A1}), which also addresses the existence and uniqueness of normalized positive minimizer $Q$ for (\ref{A1}), in the sense that $\int_{\mathbb{R}^{N}}Q^{2}dx=1$. Throughout the whole paper, for simplicity the above pair $(\mu, Q)$ is called the unique principal eigenpair of (\ref{A1}).

\begin{lem}\label{LemmaA}
For any fixed $\varepsilon>0$,  let the variational problem $\mu$ be defined by (\ref{A1}). Then we have
\begin{equation}\label{A3}
0<\mu=2\sum_{i=1}^{N}|a_i|,
\end{equation}
and
(\ref{A1}) admits a unique positive minimizer
\begin{equation}\label{0.01}
Q(x)=\Big(\prod_{i=1}^{N}\frac{2|a_{i}|}{\varepsilon\pi}\Big)^{\frac{1}{4}}
e^{-\frac{1}{\varepsilon}\sum_{j=1}^{N}|a_{j}|x_{j}^{2}}>0\ \ \text{in}\ \ \R^{N}, \ \ x=(x_1, \cdots, x_N),
\end{equation}
where  $a_i\in \R\setminus\{0\}$ is as in (\ref{1.4}) for $i=1,\cdots,N$.
\end{lem}

\noindent\textbf{Proof.}
Since the embedding $S\hookrightarrow L^{2}(\mathbb{R}^{N})$ is compact (cf. \cite[Lemma 3.1]{Z}), it is easy to check that (\ref{A1}) admits at least one minimizer $Q$ in $S$. It then follows that $\mu>0$. Additionally, we derive from \cite[Theorem 11.8]{Lieb1} that $Q$ can be chosen to be strictly positive, and $Q>0$ must be unique in $S$. Using the variational theory, we obtain that $Q$ solves the following Euler-Lagrange equation
\begin{equation}\label{3.2}
-\varepsilon\Delta Q+\frac{4}{\varepsilon}\Big(\sum_{i=1}^{N}a_{i}^{2}x_{i}^{2}\Big)Q=\mu Q\,\ \ \text{in}\, \ \mathbb{R}^{N}.
\end{equation}
It then yields from (\ref{3.2}) and the standard elliptic regularity theory that $Q\in C^{2}(\mathbb{R}^{N})$.

On the other hand, set
\begin{equation}\label{A5}
Q_{0}(x)
:=\Big(\prod_{i=1}^{N}\frac{2|a_{i}|}{\varepsilon\pi}\Big)^{\frac{1}{4}}
e^{-\frac{1}{\varepsilon}\sum_{j=1}^{N}|a_{j}|x_{j}^{2}}>0, \ \ x=(x_1, \cdots, x_N),
\end{equation}
where  $a_i\in \R\setminus\{0\}$ is as in (\ref{1.4}) for $ i=1,\cdots,N$. Note that
\begin{equation}\label{A6}
\begin{split}
\int_{\R^{N}}Q_{0}^{2}(x)dx&=\Big(\prod_{i=1}^{N}\frac{2|a_{i}|}{\varepsilon\pi}\Big)^{\frac{1}{2}}
\int_{\R^{N}}e^{-\frac{2}{\varepsilon}\sum_{j=1}^{N}|a_{j}|x_{j}^{2}}dx\\
&=\Big(\prod_{i=1}^{N}\frac{2|a_{i}|}{\varepsilon\pi}\Big)^{\frac{1}{2}}
\Big(\prod_{j=1}^{N}\int_{\R}e^{-\frac{2}{\varepsilon}|a_{j}|x_{j}^{2}}dx_{j}\Big)\\
&=\Big(\prod_{i=1}^{N}\frac{2|a_{i}|}{\varepsilon\pi}\Big)^{\frac{1}{2}}
\Big(\prod_{j=1}^{N}\sqrt{\frac{\varepsilon}{2|a_{j}|}}\int_{0}^{+\infty}
t^{-\frac{1}{2}}e^{-t}dt\Big)=1,
\end{split}
\end{equation}
and it also follows from (\ref{A1}) and (\ref{A5}) that
\begin{align}
\mu\leq&\varepsilon\int_{\R^{N}}|\nabla Q_{0}|^{2}dx+\frac{4}{\varepsilon}\int_{\R^{N}}\Big(\sum_{k=1}^{N}a_{k}^{2}x_{k}^{2}\Big)Q_{0}^{2}(x)dx\nonumber\\
=&\frac{8}{\varepsilon}\Big(\prod_{i=1}^{N}\frac{2|a_{i}|}{\varepsilon\pi}\Big)^{\frac{1}{2}}
\int_{\R^{N}}\Big(\sum_{k=1}^{N}a_{k}^{2}x_{k}^{2}\Big)
e^{-\frac{2}{\varepsilon}\sum_{j=1}^{N}|a_{j}|x_{j}^{2}}dx\nonumber\\
=&\frac{8}{\varepsilon}\Big(\prod_{i=1}^{N}\frac{2|a_{i}|}{\varepsilon\pi}\Big)^{\frac{1}{2}}
\sum_{k=1}^{N}\Big[\int_{\R}a_{k}^{2}x_{k}^{2}
e^{-\frac{2}{\varepsilon}|a_{k}|x_{k}^{2}}dx_{k}
\Big(\prod_{j\neq k}^{N}\int_{\R}
e^{-\frac{2}{\varepsilon}|a_{j}|x_{j}^{2}}dx_{j}\Big)\Big]\nonumber\\
=&\frac{8}{\varepsilon}\Big(\prod_{i=1}^{N}\frac{2|a_{i}|}{\varepsilon\pi}\Big)^{\frac{1}{2}}
\sum_{k=1}^{N}\Big[\Big(\frac{\varepsilon}{2}\Big)^{\frac{3}{2}}|a_{k}|^{\frac{1}{2}}
\Big(\int_{0}^{+\infty}
t^{\frac{1}{2}}e^{-t}dt\Big)\cdot\Big(\prod_{j\neq k}^{N}\sqrt{\frac{\varepsilon\pi}{2|a_{j}|}}\Big)\Big]\nonumber\\
=&\frac{8}{\varepsilon}\Big(\prod_{i=1}^{N}\frac{2|a_{i}|}{\varepsilon\pi}\Big)^{\frac{1}{2}}
\sum_{k=1}^{N}\Big[\frac{\varepsilon|a_{k}|}{4}\prod_{j=1}^{N}
\Big(\frac{\varepsilon\pi}{2|a_{j}|}\Big)^{\frac{1}{2}}\Big]\nonumber\\
=&2\Big(\sum_{k=1}^{N}|a_{k}|\Big).\nonumber
\end{align}
Moreover, one can check that
\begin{equation*}
-\varepsilon\Delta Q_{0}+\frac{4}{\varepsilon}\Big(\sum_{i=1}^{N}a_{i}^{2}x_{i}^{2}\Big)Q_{0}=2\Big(\sum_{i=1}^{N}|a_{i}| \Big) Q_{0}\ \,\ \text{in}\ \, \mathbb{R}^{N}.
\end{equation*}
Since (\ref{A1}) admits a unique positive minimizer $Q$, these yield that
$(\mu, Q)=\big(2\sum_{i=1}^{N}|a_{i}|, Q_{0}\big)$ is the unique principal eigenpair of the operator
$$-\varepsilon\Delta+\frac{4}{\varepsilon}\Big(\sum_{i=1}^{N}a_{i}^{2}x_{i}^{2}\Big)\,\ \ \mbox{in} \, \ L^2(\R^N),$$
which hence implies that
\begin{equation}\label{A7}
0<\mu=2\sum_{i=1}^{N}|a_{i}|,\ \ \text{and}\ \ Q(x)\equiv Q_0(x)=\Big(\prod_{i=1}^{N}\frac{2|a_{i}|}{\varepsilon\pi}\Big)^{\frac{1}{4}}
e^{-\frac{1}{\varepsilon}\sum_{j=1}^{N}|a_{j}|x_{j}^{2}}\ \ \text{in}\ \ \R^{N}.
\end{equation}
This therefore completes the proof of Lemma \ref{LemmaA}.\qed

\subsection{Proof of Theorem \ref{thm1.2}}
In this subsection, we first discuss the convergence of the unique principal eigenvalue $\lambda_{\alpha}$ for (\ref{1.1}) as $\alpha\rightarrow+\infty$, based on which we shall finally complete the proof of Theorem \ref{thm1.2}.

\begin{lem}\label{lem2.1}
Suppose $0\leq V(x)\in C^{\gamma}(\bar{\Omega})~(0<\gamma<1)$, and assume  $m(x)$ satisfies (\ref{1.4}). Then for any fixed $\varepsilon>0$, the unique principal eigenvalue $\lambda_{\alpha}$ of (\ref{1.1}) satisfies
\begin{equation}\label{2.0}
\lim_{\alpha\rightarrow+\infty}\alpha^{-1}\lambda_{\alpha}
=2\sum_{i=1}^{N}|a_{i}|>0,
\end{equation}
where $a_i\in \R\setminus\{0\}$ is as in (\ref{1.4}) for $ i=1,\cdots,N$.
\end{lem}

\noindent\textbf{Proof.}
Under the assumptions of Lemma \ref{lem2.1}, suppose $(\lambda_\alpha, u_\alpha)$ is the unique principal eigenpair of (\ref{1.1}). We then deduce from (\ref{1.02}) and (\ref{A1}) that
\begin{equation}\label{2.1}
\begin{split}
\lambda_{\alpha}
=&\int_{\Omega}\varepsilon|\nabla u_{\alpha}|^2dx+\int_{\Omega}\Big[\frac{\alpha^{2}}{\varepsilon}|\nabla m|^{2}
+\alpha \Delta m
+V(x)\Big]u_{\alpha}^2dx\\
=&\int_{\Omega}\varepsilon|\nabla u_{\alpha}|^2dx+\int_{\Omega}\Big[\frac{4\alpha^{2}}{\varepsilon}
\Big(\sum_{i=1}^{N}a_{i}^{2}x_{i}^{2}\Big)+V(x)\Big]u_{\alpha}^2dx\\
\geq&\alpha\Big\{\varepsilon\int_{\R^{N}}|\nabla \tilde{v}_{\alpha}|^{2}dx
+\frac{4}{\varepsilon}\int_{\R^{N}}\Big(\sum_{i=1}^{N}a_{i}^{2}x_{i}^{2}\Big)
\tilde{v}_{\alpha}^{2}dx\Big\}\\
\geq&\alpha\mu\ \ \text{as}\ \ \alpha\rightarrow+\infty,
\end{split}
\end{equation}
where $\mu=2\sum_{i=1}^{N}|a_{i}|>0$ is as in (\ref{A3}),  $\tilde{v}_{\alpha}(x)=\alpha^{-\frac{N}{4}}u_{\alpha}(\alpha^{-\frac{1}{2}}x)>0$, and
$u_{\alpha}(x)\equiv 0$ in $\R^{N}\backslash\Omega$. This proves the lower bound of $\lambda_{\alpha}$ as $\alpha\rightarrow+\infty$.

We next analyze the upper estimate of $\lambda_{\alpha}$ as $\alpha\rightarrow+\infty$. Choose a non-negative cut off function $\xi(x)\in C_{0}^{\infty}(\R^{N})$, where $\xi(x)=1$ for $|x|\leq r_{1}$ and $\xi(x)=0$ for $|x|\geq 2r_{1}$, where $r_{1}>0$ is so small  that $B_{2r_{1}}(0)\subset\Omega$. For any $\alpha>0$, define
\begin{equation}\label{2.2}
\hat{\phi}_{\alpha}(x):=\hat{A}_{\alpha}\alpha^{\frac{N}{4}}
\xi(x)Q(\alpha^{\frac{1}{2}}x),\ \ x\in\Omega,
\end{equation}
where $Q(x)>0$ is as in (\ref{0.01}), and $\hat{A}_{\alpha}>0$ is chosen such that $\int_{\Omega}|\hat{\phi}_{\alpha}|^{2}dx=1$. By the exponential decay of $Q(x)$, we then have
\begin{equation}\label{2.3}
\frac{1}{\hat{A}_{\alpha}^{2}}=1+O(\alpha^{-\infty})
\ \ \text{as}\ \ \alpha\rightarrow+\infty.
\end{equation}
Here and below we use $f(t)=O(t^{-\infty})$ to denote a function $f$ satisfying $\lim_{t\rightarrow+\infty}|f(t)|t^{s}=0$ for all $s>0$.
Applying Lemma \ref{LemmaA}, one can calculate from (\ref{1.4}), (\ref{1.02}), (\ref{2.2}) and (\ref{2.3}) that
\begin{equation}\label{2.4}
\begin{split}
\lambda_{\alpha}
\leq&\int_{\Omega}\varepsilon|\nabla \hat{\phi}_{\alpha}|^2dx+\int_{\Omega}\Big[\frac{\alpha^{2}}{\varepsilon}|\nabla m|^{2}
+\alpha \Delta m
+V(x)\Big]\hat{\phi}_{\alpha}^2dx\\
=&\int_{\Omega}\varepsilon|\nabla \hat{\phi}_{\alpha}|^2dx
+\int_{\Omega}\Big[\frac{4\alpha^{2}}{\varepsilon}
\Big(\sum_{i=1}^{N}a_{i}^{2}x_{i}^{2}\Big)+V(x)\Big]\hat{\phi}_{\alpha}^2dx\\
=&\alpha\Big[\int_{\R^{N}}\varepsilon|\nabla Q|^{2}dx+\frac{4}{\varepsilon}\int_{\R^{N}}\Big(\sum_{i=1}^{N}a_{i}^{2}x_{i}^{2}\Big)Q^{2}(x)dx+o(1)\Big]\\
=&\alpha\big[\mu+o(1)\big]
\ \ \text{as}\ \ \alpha\rightarrow+\infty,
\end{split}
\end{equation}
where $\mu=2\sum_{i=1}^{N}|a_{i}|>0$ is as in (\ref{A3}).
We thus conclude from (\ref{A3}), (\ref{2.1}) and (\ref{2.4}) that (\ref{2.0}) holds true, and the proof of Lemma \ref{lem2.1} is therefore complete.\qed

Applying Lemmas \ref{LemmaA} and \ref{lem2.1}, we next analyze the convergence of the unique normalized principal eigenfunction $u_{\alpha}$ for (\ref{1.1}) as $\alpha\rightarrow+\infty$, which then completes the proof of Theorem \ref{thm1.2}.

\begin{prop}\label{pro2.2}
Suppose $0\leq V(x)\in C^{\gamma}(\bar{\Omega})~(0<\gamma<1)$, and assume  $m(x)$ satisfies (\ref{1.4}). Then for any fixed $\varepsilon>0$, the unique normalized principal eigenfunction $u_{\alpha}$ of (\ref{1.1}) satisfies
\begin{equation}\label{2.5m}
\begin{split}
\hat{w}_\alpha(x)&:=\alpha^{-\frac{N}{4}}u_{\alpha}\big(\alpha^{-\frac{1}{2}}x+d_\alpha\big)\rightarrow Q(x)=\Big(\prod_{i=1}^{N}\frac{2|a_{i}|}{\varepsilon\pi}\Big)^{\frac{1}{4}}
e^{-\frac{1}{\varepsilon}\sum_{j=1}^{N}|a_{j}|x_{j}^{2}}>0\\
& \ \, \text{strongly in}\ \, H^{1}(\mathbb{R}^{N})\cap L^{\infty}(\mathbb{R}^{N})
\ \, \text{as}\ \, \alpha\rightarrow+\infty,
\end{split}
\end{equation}
where $a_i\in \R\setminus\{0\}$ is as in (\ref{1.4}) for $ i=1,\cdots,N$, $u_{\alpha}(x)\equiv0$ in $\R^{N}\backslash\Omega$,   and $d_\alpha\in \Omega$ is a global maximum point of $u_\alpha$ and satisfies
\begin{equation}\label{2.0021}
\lim_{\alpha\rightarrow+\infty}\alpha^{\frac{1}{2}}|d_\alpha|=0.
\end{equation}
\end{prop}

\noindent\textbf{Proof.}
For any $\alpha>0$, set
\begin{equation}\label{2.19}
\hat{w}_\alpha(x):=\left\{
\begin{split}\alpha^{-\frac{N}{4}}u_{\alpha}(\alpha^{-\frac{1}{2}}x+d_\alpha),\ \ &x\in\hat{\Omega}_{\alpha}:=\{x\in\R^{N}:\ \alpha^{-\frac{1}{2}}x+d_\alpha\in\Omega\};\\[2mm]
0, \qquad\qquad &x\in\mathbb{R}^{N}\backslash\hat{\Omega}_{\alpha},
\end{split}
\right.
\end{equation}
where $d_\alpha\in\Omega$ is a global maximum point of $u_\alpha$. It is easy to check that
\begin{equation}\label{2.22}
\int_{\R^{N}}\hat{w}_{\alpha}^2(x)dx=\int_{\Omega}u_{\alpha}^2(x)dx=1\ \ \text{for any}\ \, \alpha>0.
\end{equation}
One can also calculate from (\ref{1.4}), (\ref{1.31}) and (\ref{2.19}) that $\hat{w}_{\alpha}$ satisfies
\begin{equation}\label{2.20}
-\varepsilon\Delta \hat{w}_\alpha+\Big[\frac{4}{\varepsilon}\sum_{i=1}^{N}a_{i}^{2}\big(x_{i}+\alpha^{\frac{1}{2}}
(d_\alpha)_{i}\big)^{2}
+\alpha^{-1}V\big(\alpha^{-\frac{1}{2}}x+d_\alpha\big)\Big]\hat{w}_\alpha
=\alpha^{-1}\lambda_{\alpha}\hat{w}_\alpha\ \ \text{in}\ \ \hat{\Omega}_{\alpha},
\end{equation}
where $x=(x_{1},\cdots,x_{N})\in\hat{\Omega}_{\alpha}$, and $d_\alpha=((d_\alpha)_{1},\cdots,(d_\alpha)_{N})\in\Omega$. Since the origin is always a global maximum point of $\hat{w}_{\alpha}(x)$, we have
$$-\Delta \hat{w}_{\alpha}(0)\geq0\ \ \text{and}\ \ \hat{w}_{\alpha}(0)>0\ \
\text{for any}\ \, \alpha>0,$$
which implies from (\ref{2.20}) that
\begin{equation}\label{2.20m}
0\leq\frac{4\alpha}{\varepsilon}\sum_{i=1}^{N}a_{i}^{2}(d_\alpha)_{i}^{2}\leq\alpha^{-1}\lambda_{\alpha}\ \
\text{for any}\ \, \alpha>0.
\end{equation}
One can further obtain from (\ref{1.4}), (\ref{2.0}) and (\ref{2.20m}) that there exist a sequence $\{\alpha_n\}$, which satisfies $\lim\limits_{n\rightarrow\infty}\alpha_n=+\infty$, and a point $x_0\in \R^{N}$ such that
\begin{equation}\label{2.21}
\lim_{n\rightarrow\infty}\alpha_{n}^{\frac{1}{2}}d_{\alpha_n}=x_0.
\end{equation}
It then yields from (\ref{2.21}) that $d_{\alpha_n}$ converges to the origin as $n\rightarrow\infty$. Since the origin is an interior point of $\Omega$,   the definition of $\hat{\Omega}_{\alpha_n}$ in (\ref{2.19}) yields that $\lim\limits_{n\rightarrow\infty}\hat{\Omega}_{\alpha_n}=\R^{N}$.

We now derive from (\ref{2.0}) and (\ref{2.19})--(\ref{2.21}) that
\begin{equation}\label{2.23}
0\leq \varepsilon\int_{\R^N}|\nabla \hat{w}_{\alpha_n}|^2dx\leq2\sum_{i=1}^{N}|a_{i}|+1\ \ \text{as}\ \ n\rightarrow\infty,
\end{equation}
and
\begin{align}
0\leq& \int_{\R^{N}}
\Big(\sum_{i=1}^{N}a_{i}^{2}x_{i}^{2}\Big)\hat{w}_{\alpha_n}^2(x)dx\nonumber\\
\leq&2\sum_{i=1}^{N}a_{i}^{2}\int_{\R^{N}}
\big(x_{i}+\alpha_n^{\frac{1}{2}}(d_{\alpha_n})_{i}\big)^{2}
\hat{w}_{\alpha_n}^2(x)dx
+2\alpha_n\sum_{i=1}^{N}a_{i}^{2}\int_{\R^N}(d_{\alpha_n})_{i}^{2}\hat{w}_{\alpha_n}^2(x)dx\label{2.24}\\
\leq& C\Big(\sum_{i=1}^{N}|a_{i}|\Big)
\ \ \text{as}\ \ \alpha_n\rightarrow\infty,\nonumber
\end{align}
where  $a_i\in \R\setminus\{0\}$ is as in (\ref{1.4}) for $ i=1,\cdots,N$, and $d_{\alpha_n}=((d_{\alpha_n})_{1},\cdots,(d_{\alpha_n})_{N})\in\Omega$.
It then follows from (\ref{2.22}), (\ref{2.23}) and (\ref{2.24}) that $\{\hat{w}_{\alpha_n}\}$ is bounded uniformly in $S$, which is defined by (\ref{A2}). Since the embedding $S\hookrightarrow L^{q}(\R^{N})$ is compact (cf. \cite[Lemma 3.1]{Z}) for $2\leq q<2^{*}$, where $2^{*}=+\infty$ if $N=1,2$ and $2^{*}=\frac{2N}{N-2}$ if $N\geq3$, we obtain that there exist a subsequence, still denoted by $\{\hat{w}_{\alpha_n}\}$, of $\{\hat{w}_{\alpha_n}\}$ and a function $0\leq\hat{w}_{0}\in S$ such that
\begin{equation}\label{2.25}
\begin{split}
\hat{w}_{\alpha_{n}}\rightharpoonup\hat{w}_{0}\ \ \text{weakly in}\ \, H^{1}(\mathbb{R}^{N})\ \ \text{as}\ \ n\rightarrow\infty,\qquad\qquad\\
\hat{w}_{\alpha_{n}}\rightarrow\hat{w}_{0} \ \ \text{strongly in}\ \, L^{q}(\mathbb{R}^{N})\ \ \text{as}\ \ n\rightarrow\infty,\ \  2\leq q<2^*.
\end{split}
\end{equation}
We thus obtain from (\ref{2.22}) and (\ref{2.25}) that $\int_{\R^N}|\hat{w}_0|^2dx=1$.

Applying (\ref{A3}) and (\ref{2.0}), together with Fatou's Lemma, we derive from (\ref{A1}), (\ref{2.22}), (\ref{2.20}), (\ref{2.21}) and (\ref{2.25}) that
\begin{align}
2\sum_{i=1}^{N}|a_{i}|=&\lim_{n\rightarrow\infty}\alpha_{n}^{-1}\lambda_{\alpha_n}\int_{\R^N}\hat{w}_{\alpha_{n}}^2dx\nonumber\\
=&\lim_{n\rightarrow\infty}\Big\{\int_{\R^N}\varepsilon|\nabla \hat{w}_{\alpha_n}|^2dx+
\int_{\R^N}\Big[\frac{4}{\varepsilon}\sum_{i=1}^{N}a_{i}^{2}\big(x_{i}+{\alpha}_n^{\frac{1}{2}}
(d_{\alpha_n})_{i}\big)^{2}\nonumber\\
&\qquad\quad+\alpha_n^{-1}V({\alpha}_n^{-\frac{1}{2}}x+d_{\alpha_n})\Big]
\hat{w}_{\alpha_n}^2dx\Big\}\label{2.9m}\\
\geq&\liminf_{n\rightarrow\infty}\int_{\R^N}\varepsilon|\nabla \hat{w}_{\alpha_n}|^2dx
+\liminf_{n\rightarrow\infty}\frac{4}{\varepsilon}\int_{\R^N}\sum_{i=1}^{N}a_{i}^{2}
\big(x_{i}+{\alpha}_n^{\frac{1}{2}}
(d_{\alpha_n})_{i}\big)^{2}\hat{w}_{\alpha_n}^2dx\nonumber\\
\geq&\varepsilon\int_{\R^N}|\nabla \hat{w}_0|^2dx+\frac{4}{\varepsilon}\int_{\R^N}\sum_{i=1}^{N}a_{i}^{2}
\big(x_{i}+(x_{0})_i\big)^{2}\hat{w}_{0}^2(x)dx\nonumber\\
=&\varepsilon\int_{\R^N}|\nabla \hat{w}_0|^2dx+\frac{4}{\varepsilon}\int_{\R^N}\Big(\sum_{i=1}^{N}a_{i}^{2}
x_i^{2}\Big)\hat{w}_{0}^2(x-x_0)dx
\geq\mu=2\sum_{i=1}^{N}|a_{i}|,\nonumber
\end{align}
where $x_0=((x_0)_{1},\cdots,(x_0)_{N})\in\R^N$ is as in (\ref{2.21}), and $\hat{w}_0(x)\geq0$ is given by (\ref{2.25}). We thus obtain from Lemma \ref{LemmaA} that $\hat{w}_0(x-x_0)\equiv Q(x)>0$ in $\R^N$, where $Q(x)$ is as in (\ref{0.01}). Moreover, since the definition of $\hat{w}_{\alpha_n}(x)$ in (\ref{2.19}) implies that the origin is always a global maximum point of $\hat{w}_{\alpha_n}(x)$ for all $n>0$,  it  follows that the origin is also a global maximum point of $\hat{w}_{0}(x)$ in $\R^N$.
On the other hand, it yields from (\ref{0.01}) that the origin is the unique global maximum point of $Q(x)$. We thus conclude from above   that $|x_0|=0$, and hence (\ref{2.0021}) holds true in view of (\ref{2.21}).

We next claim that
\begin{equation}\label{2.10}
\hat{w}_{\alpha_{n}}(x)\rightarrow \hat{w}_{0}(x)\equiv Q(x)\ \ \text{strongly in}\ \ H^{1}(\mathbb{R}^{N})\ \ \text{as}\ \ n\rightarrow\infty.
\end{equation}
Indeed, it yields from (\ref{2.0}) and (\ref{2.20}) that for any fixed $\varepsilon>0$,
\begin{equation}\label{2.11}
-\Delta \hat{w}_{\alpha_{n}}-\frac{4}{\varepsilon}\Big(\sum_{i=1}^{N}|a_{i}|\Big)\hat{w}_{\alpha_{n}}\leq0\ \ \text{in}\ \ \hat{\Omega}_{\alpha_{n}}\ \ \text{as}\ \
n\rightarrow\infty,
\end{equation}
where $a_i\in \R\setminus\{0\}$ is as in (\ref{1.4}) for $ i=1,\cdots,N$.
Since $\lim\limits_{n\rightarrow\infty}\hat{\Omega}_{\alpha_{n}}=\R^{N}$, we   have $B_{R+2}(0)\subseteq\hat{\Omega}_{\alpha_{n}}$ as $n\rightarrow\infty$ for any $R>0$. By De Giorgi-Nash-Moser theory  (cf. \cite[Theorem 4.1]{Han}), we then obtain from (\ref{2.11}) that there exists a constant $C>0$, independent of $n>0$ and above $R>0$, such that for any $z\in B_{R+1}(0)$,
\begin{equation}\label{2.12}
\max_{B_{\frac{1}{2}}(z)}\hat{w}_{\alpha_{n}}(x)\leq C\Big(\int_{{B_{1}(z)}}|\hat{w}_{\alpha_{n}}(x)|^{2} dx\Big)^\frac{1}{{2}}\ \ \mbox{as}\ \ n\rightarrow\infty.
\end{equation}
We thus derive from (\ref{2.22}) and (\ref{2.12}) that
\begin{equation}\label{2.13}
\{\hat{w}_{\alpha_{n}}\}\ \ \text{is bounded uniformly in}\, \ H^{1}(B_{R+1}(0))\cap L^{\infty}(B_{R+1}(0)) \ \ \mbox{as}\ \ n\rightarrow\infty.
\end{equation}

Moreover, it follows from (\ref{1.4}), (\ref{2.0}) and (\ref{2.20}) that there exists a large constant $R_{1}>0$, independent of $n>0$, such that
\begin{equation}\label{2.14}
-\Delta \hat{w}_{\alpha_{n}}+\hat{w}_{\alpha_{n}}\leq0\ \ \text{uniformly in}\ \ \hat{\Omega}_{\alpha_{n}}\backslash B_{R_{1}}(0) \ \ \mbox{as}\ \ n\rightarrow\infty.
\end{equation}
Applying the comparison principle, one can obtain from (\ref{2.13}) and (\ref{2.14}) that there exists a constant $C>0$, independent of $n>0$, such that
\begin{equation}\label{2.15}
|\hat{w}_{\alpha_{n}}(x)|\leq Ce^{-|x|}\ \ \text{uniformly in}\ \ \hat{\Omega}_{\alpha_{n}}\backslash B_{R_{1}}(0) \ \ \mbox{as}\ \ n\rightarrow\infty.
\end{equation}
It then follows from (\ref{2.0021}), (\ref{2.25}), (\ref{2.13}) and (\ref{2.15}) that
\begin{equation}\label{m1}
\lim_{n\rightarrow\infty}\int_{\R^N}\sum_{i=1}^{N}a_{i}^{2}
\big(x_{i}+{\alpha}_n^{\frac{1}{2}}
(d_{\alpha_n})_{i}\big)^{2}\hat{w}_{\alpha_n}^2dx
=\int_{\R^{N}}\Big(\sum_{i=1}^{N}a_{i}^{2}x_{i}^{2}\Big)\hat{w}_{0}^{2}(x)dx,
\end{equation}
and
\begin{equation}\label{m2}
\lim_{n\rightarrow\infty}\int_{\R^{N}}\alpha_n^{-1}V\big({\alpha}_n^{-\frac{1}{2}}x+d_{\alpha_n}\big)
\hat{w}_{\alpha_n}^2dx=0.
\end{equation}
We thus obtain from (\ref{2.9m}), (\ref{m1}) and (\ref{m2}) that
$$\lim_{n\rightarrow\infty}\int_{\R^{N}}|\nabla \hat{w}_{\alpha_{n}}|^{2}dx
=\int_{\R^{N}}|\nabla \hat{w}_{0}|^{2}dx=\int_{\R^{N}}|\nabla Q|^{2}dx,$$
which further proves the claim (\ref{2.10}) in view of (\ref{2.25}).

We finally claim that
\begin{equation}\label{2.16}
\hat{w}_{\alpha_{n}}(x)\rightarrow Q(x)\ \ \text{strongly in}\ \ L^{\infty}(\R^{N})\ \ \text{as}\ \ n\rightarrow\infty.
\end{equation}
Actually, (\ref{2.16}) holds true for $N=1$ in view of (\ref{2.10}), together with the fact that $H^{1}(\R)\hookrightarrow L^{\infty}(\R)$ is continuous (cf. \cite[Theorem 8.5]{Lieb1}). As for the case $N\geq2$, we just need to show that \begin{equation}\label{2.17}
\hat{w}_{\alpha_{n}}(x)\rightarrow Q(x)\ \ \text{strongly in}\ \ L_{loc}^{\infty}(\R^{N})\ \ \text{as}\ \ n\rightarrow\infty
\end{equation}
in view of (\ref{0.01}) and (\ref{2.15}). Indeed, defining
$$F_{\alpha_n}(x):=\Big[\frac{\lambda_{\alpha_n}}{\varepsilon\alpha_n}
-\frac{4}{\varepsilon^{2}}\Big(\sum_{i=1}^{N}a_{i}^{2}\big(x_{i}+\alpha_{n}^{\frac{1}{2}}
(d_{\alpha_n})_{i}\big)^{2}\Big)
-\frac{1}{\varepsilon\alpha_n}V(\alpha_n^{-\frac{1}{2}}x+d_{\alpha_n})\Big]\hat{w}_{\alpha_n}(x),$$
it follows from (\ref{2.20}) that
\begin{equation}\label{2.18}
-\Delta \hat{w}_{\alpha_n}=F_{\alpha_n}(x)\ \ \text{in}\ \ \hat{\Omega}_{\alpha_n},
\end{equation}
where the domain $\hat{\Omega}_{\alpha_n}$ is as in (\ref{2.19}). By the standard elliptic regularity theory, we then deduce that there exists a subsequence, still denoted by $\{\hat{w}_{\alpha_n}\}$, of $\{\hat{w}_{\alpha_n}\}$ such that (\ref{2.17}) holds true. This hence proves the claim (\ref{2.16}).

Because the positive function $Q(x)>0$ is unique in view of Lemma \ref{LemmaA} and the above results are independent of the subsequence that we choose, we conclude that the above results hold true for the whole sequence $\{\hat{w}_{\alpha}\}$. This therefore completes the proof of Proposition \ref{pro2.2} in view of (\ref{0.01}).
\qed

Theorem \ref{thm1.2} now follows directly from Lemma  \ref{lem2.1} and Proposition \ref{pro2.2}.
Before ending this section, we discuss the uniqueness of the global maximum point $d_\alpha\in \Omega$ for   $u_\alpha$ as $\alpha\rightarrow+\infty$, for the case where $|a_i|=|a_j|>0$ in (\ref{1.4}) holds for all $i\neq j$.

\begin{cor}\label{cor2.3}
Under the assumptions of Proposition \ref{pro2.2}, suppose $|a_i|=|a_j|>0$ in (\ref{1.4}) holds for all $i\neq j$. Then the global maximum point $d_\alpha\in \Omega$ of $u_\alpha$ for (\ref{1.1}) must be unique as $\alpha\rightarrow+\infty$.
\end{cor}

\noindent\textbf{Proof.}
Since $|a_i|=|a_j|>0$ in (\ref{1.4}) holds for all $i\neq j$, we obtain from Proposition \ref{pro2.2} that
\begin{equation}\label{c1}
\begin{split}
\hat{w}_\alpha(x)&:=\alpha^{-\frac{N}{4}}u_{\alpha}(\alpha^{-\frac{1}{2}}x+d_\alpha)
\rightarrow \Big(\frac{2|a_{1}|}{\varepsilon\pi}\Big)^{\frac{N}{4}}
e^{-\frac{|a_{1}|}{\varepsilon}|x|^{2}}>0\\
&\ \ \text{strongly in}\ \ H^{1}(\R^{N})\cap L^{\infty}(\R^{N})\ \ \text{as}\ \ \alpha\rightarrow+\infty.
\end{split}
\end{equation}
Moreover, by the standard elliptic regularity theory, one can derive from (\ref{2.18}) that there exists a subsequence $\{\hat{w}_{\alpha_l}\}$ of $\{\hat{w}_\alpha\}$ such that for some $\hat{Q}\in C_{loc}^2(\R^N)$,
\begin{equation}\label{c2}
\hat{w}_{\alpha_l}(x)
\rightarrow \hat{Q}(x) \ \ \text{in}\ \ C_{loc}^2(\R^N)\ \ \text{as}\ \ \alpha_l\rightarrow\infty.
\end{equation}
It then yields from (\ref{c1}) and (\ref{c2}) that
\begin{equation}\label{c3}
\hat{w}_{\alpha_l}(x)
\rightarrow \hat{Q}(x)\equiv\Big(\frac{2|a_{1}|}{\varepsilon\pi}\Big)^{\frac{N}{4}}
e^{-\frac{|a_{1}|}{\varepsilon}|x|^{2}}>0 \ \ \text{in}\ \ C_{loc}^2(\R^N)\ \ \text{as}\ \ \alpha_l\rightarrow\infty.
\end{equation}

One can check that there exists a small constant $\hat{\delta}>0$ such that $\hat{Q}''(r)<0$  holds for all $0\leq r<\hat{\delta}$.
Furthermore, since the origin is the unique global maximum point of $\hat{Q}(x)$, we derive from (\ref{2.15}) and (\ref{c3}) that all global maximum points of $\hat{w}_{\alpha_l}$ must stay in a small ball $B_{\hat{\delta}}(0)$ as $\alpha_l\rightarrow\infty$.
We then deduce from \cite[Lemma 4.2]{Ni} that $\hat{w}_{\alpha_l}$ has a unique maximum point for sufficiently large $\alpha_l>0$, which is just the origin. This proves the uniqueness of global maximum points for $\hat{w}_{\alpha_l}$ as $\alpha_l\rightarrow\infty$.
Because the above uniqueness is independent of the subsequence $\{\hat{w}_{\alpha_{l}}\}$ that we choose, we conclude that it holds essentially true for the whole sequence $\{\hat{w}_{\alpha}\}$. The proof of Corollary \ref{cor2.3} is therefore complete.
\qed

\section{Refined Limiting Profiles of $(\lambda_{\alpha}, u_{\alpha})$ as $\alpha\rightarrow+\infty$}

In this section, without special notations we shall always assume that $V(x)$  satisfies $V(x)\in C^{2}(B_{r_0}(0))$ in some small ball $B_{r_0}(0)\subset\Omega$, where $r_0>0$. Following the convergence of Theorem \ref{thm1.2}, the main purpose of this section is to prove Theorem \ref{thm1.3} on the   refined limiting profiles of the unique principal eigenpair $(\lambda_{\alpha}, u_{\alpha})$ for (\ref{1.1}) as $\alpha\rightarrow+\infty$.

We first give some further remarks on the unique normalized principal eigenfunction $u_{\alpha}>0$ of (\ref{1.1}).
For any $\alpha>0$, set
\begin{equation}\label{2.6}
w_\alpha(x):=\left\{
\begin{split}\alpha^{-\frac{N}{4}}u_{\alpha}(\alpha^{-\frac{1}{2}}x),\ \ &x\in\Omega_{\alpha}:=\big\{x\in\R^{N}:\ \alpha^{-\frac{1}{2}}x\in\Omega\big\};\\[2mm]
0, \qquad\qquad &x\in\mathbb{R}^{N}\backslash\Omega_{\alpha}.
\end{split}
\right.
\end{equation}
It is easy to check that
\begin{equation}\label{M2.7}
\int_{\R^{N}}w_\alpha^{2}(x)dx=\int_{\Omega}u_\alpha^{2}(x)dx=1\ \ \text{for all}\ \ \alpha>0.
\end{equation}
One can also calculate from (\ref{1.4}), (\ref{1.31}) and (\ref{2.6}) that $w_\alpha(x)$ satisfies
\begin{equation}\label{2.7}
-\varepsilon\Delta w_\alpha+\Big[\frac{4}{\varepsilon}\Big(\sum_{i=1}^{N}a_{i}^{2}x_{i}^{2}\Big)
+\alpha^{-1}V(\alpha^{-\frac{1}{2}}x)\Big]w_\alpha=\alpha^{-1}\lambda_{\alpha}w_\alpha\ \,\ \text{in}\ \, \Omega_{\alpha},
\end{equation}
where  the domain $\Omega_{\alpha}$ is as in (\ref{2.6}).
Since the origin is an interior point of $\Omega$,  the definition of $\Omega_{\alpha}$ in (\ref{2.6}) gives that $\lim\limits_{\alpha\rightarrow+\infty}\Omega_{\alpha}=\R^{N}$.
By the comparison principle, similar to (\ref{2.15}), one can derive  that there exists a constant $C>0$, independent of $\alpha>0$, such that for any fixed large $R_1>0,$
\begin{equation}\label{2.15m}
|w_{\alpha}(x)|\leq Ce^{-|x|}\ \ \text{uniformly in}\ \ \Omega_{\alpha}\backslash B_{R_{1}}(0) \ \ \mbox{as}\ \ \alpha\rightarrow+\infty.
\end{equation}
Additionally, we obtain from Theorem \ref{thm1.2} that
\begin{equation}\label{2.5}
\begin{split}
w_{\alpha}(x)&\rightarrow Q(x)=\Big(\prod_{i=1}^{N}\frac{2|a_{i}|}{\varepsilon\pi}\Big)^{\frac{1}{4}}
e^{-\frac{1}{\varepsilon}\sum_{j=1}^{N}|a_{j}|x_{j}^{2}}>0 \\
 \ \ &\text{strongly in}\ L^{\infty}(\mathbb{R}^{N})
\ \ \text{as}\ \ \alpha\rightarrow+\infty,
\end{split}
\end{equation}
where $w_{\alpha}(x)$ is defined by (\ref{2.6}), and $a_{i}\in \R\setminus\{0\}$ is as in (\ref{1.4}) for $i=1,\cdots,N$.

Under a weaker regularity assumption on $V(x)$, we now analyze the second-order asymptotic expansion of $\lambda_{\alpha}$ as $\alpha\rightarrow+\infty$, for which we shall establish an important Pohozaev type identity.

\begin{lem}\label{lem3.1}
Suppose $m(x)$ satisfies (\ref{1.4}), and assume $0\leq V (x)\in C^{\gamma}(\bar{\Omega})~(0<\gamma<1)$. Then for any fixed $\varepsilon>0$, the unique principal eigenvalue $\lambda_{\alpha}$ of (\ref{1.1}) satisfies
\begin{equation}\label{3.0}
\lambda_{\alpha}=2\alpha\sum_{i=1}^{N}|a_{i}|+V(0)+o(1)\ \ \text{as}\ \ \alpha\rightarrow+\infty,
\end{equation}
where $a_{i}\in \R\setminus\{0\}$ is as in (\ref{1.4}) for $i=1,\cdots,N$.
\end{lem}

\noindent\textbf{Proof.}
Define
\begin{equation}\label{3.1}
\mathcal{N}:=-\varepsilon\Delta+\frac{4}{\varepsilon}
\sum_{i=1}^{N}a_{i}^{2}x_{i}^{2}-2\sum_{i=1}^{N}|a_{i}|\ \ \text{in}\ \ \R^{N},\ \ x=(x_{1},\cdots,x_{N}),
\end{equation}
where $\varepsilon>0$ is fixed, and $a_{i}\in \R\setminus\{0\}$ is as in (\ref{1.4}) for $i=1,\cdots,N$.
We first claim that
\begin{equation}\label{3.4}
\ker\mathcal{N}=span\big\{Q\big\},
\end{equation}
where $Q(x)>0$ is as in (\ref{0.01}).
Indeed, it follows from (\ref{A3}), (\ref{3.2}) and (\ref{3.1}) that $Q(x)$ satisfies
\begin{equation}\label{3.3}
\mathcal{N}Q=0\ \,\ \text{in}\, \ \R^{N}.
\end{equation}
On the other hand, assume that there exists a nonzero function $u_0(x)\in\ker\mathcal{N}$, $i.e.$, $\mathcal{N}u_0=0$ in $\R^N$. We then derive from (\ref{A3}) and (\ref{3.1}) that
$$\mu=2\sum_{i=1}^{N}|a_{i}|=\frac{1}{\|u_0\|^2_{L^2(\R^N)}}\Big\{\int_{\R^N}
\varepsilon|\nabla u_0|^2dx+\frac{4}{\varepsilon}
\int_{\R^N}\Big(\sum_{i=1}^{N}a_{i}^{2}x_{i}^{2}\Big) u_0^2dx\Big\},$$
which thus implies from Lemma \ref{LemmaA} that there exists a constant $c_0\in\R$ such that $u_0=c_0Q$ in $\R^N$.
We hence conclude from (\ref{3.3}) that the claim (\ref{3.4}) holds true.

Following (\ref{2.5}), we next set
\begin{equation}\label{3.5}
z_{\alpha}(x):=w_{\alpha}(x)-Q(x)\ \ \text{in}\ \ \R^N,
\end{equation}
where $w_{\alpha}(x)$ is as in (\ref{2.6}), and $Q(x)>0$ is as in (\ref{0.01}). It then follows from (\ref{2.7}), (\ref{3.1}), (\ref{3.3}) and (\ref{3.5}) that $z_{\alpha}(x)$ satisfies
\begin{equation}\label{3.6}
\mathcal{N}z_{\alpha}=\Big[-\alpha^{-1}V(\alpha^{-\frac{1}{2}}x)
+\Big(\alpha^{-1}\lambda_\alpha-2\sum_{i=1}^{N}|a_{i}|\Big)\Big]w_\alpha\ \,\ \text{in}\ \, \Omega_\alpha,
\end{equation}
where $a_{i}\in\R\setminus\{0\}$ is as in (\ref{1.4}) for $i=1,\cdots,N$, and the domain $\Omega_\alpha$ is as in (\ref{2.6}). Moreover, applying the local elliptic estimate (cf. \cite[(3.15)]{GT}) to (\ref{2.7}), we derive from (\ref{2.0}) and (\ref{2.15m}) that there exist constants $C>0$, independent of $\alpha>0$, such that for sufficiently large $R>0,$
\begin{equation}\label{3.7}
|\nabla w_\alpha(x)|\leq Ce^{-\frac{|x|}{2}}\ \ \text{uniformly in}\ \ \Omega_\alpha\backslash B_R(0)\ \ \text{as}\ \ \alpha\rightarrow+\infty.
\end{equation}

Since $\int_{\R^N}Q^2(x)dx=1$, we now calculate from (\ref{0.01}), (\ref{2.15m}), (\ref{2.5}), (\ref{3.3}) and (\ref{3.7}) that
\begin{equation}\label{3.8}
\begin{split}
&\int_{B_{R_{0}\alpha^{\frac{1}{4}}}(0)}(\mathcal{N}z_{\alpha})Qdx\\
=&\int_{B_{R_{0}\alpha^{\frac{1}{4}}}(0)}(\mathcal{N}Q)z_{\alpha}dx
+\varepsilon\int_{\partial B_{R_{0}\alpha^{\frac{1}{4}}}(0)}\Big(z_{\alpha}\frac{\partial Q}{\partial \nu}-
Q\frac{\partial z_\alpha}{\partial \nu}\Big)dS\\
=&o\big(e^{-R_{0}\alpha^{\frac{1}{4}}}\big)\ \ \text{as}\ \ \alpha\rightarrow+\infty,
\end{split}
\end{equation}
and
\begin{equation}\label{3.9}
\begin{split}
&\int_{B_{R_{0}\alpha^{\frac{1}{4}}}(0)}\Big[-\alpha^{-1}V(\alpha^{-\frac{1}{2}}x)
+\Big(\alpha^{-1}\lambda_{\alpha}-2\sum_{i=1}^{N}|a_{i}|\Big)\Big]w_\alpha Qdx\\
=&-\alpha^{-1}\big[V(0)+o(1)\big]+\Big(\alpha^{-1}\lambda_{\alpha}-2\sum_{i=1}^{N}|a_{i}|\Big)\big[1+o(1)\big]\ \ \text{as}\ \ \alpha\rightarrow+\infty,
\end{split}
\end{equation}
where $R_0>0$ is small enough, and $\nu=(\nu^{1},\cdots,\nu^{N})$ denotes the outward unit normal vector of $\partial B_{R_{0}\alpha^{\frac{1}{4}}}(0)$. We thus conclude from (\ref{3.6}) and (\ref{3.8}) that
\begin{equation}\label{m}
-\alpha^{-1}\big[V(0)+o(1)\big]+\Big(\alpha^{-1}\lambda_{\alpha}-2\sum_{i=1}^{N}|a_{i}|\Big)[1+o(1)]
=o(e^{-R_{0}\alpha^{\frac{1}{4}}})\ \ \text{as}\ \ \alpha\rightarrow+\infty.
\end{equation}
It further follows from (\ref{m}) that
\begin{equation}\label{m11}
\alpha^{-1}\lambda_{\alpha}-2\sum_{i=1}^{N}|a_{i}|=\alpha^{-1}\big[V(0)+o(1)\big]\ \ \text{as}\ \ \alpha\rightarrow+\infty,
\end{equation}
which hence implies that (\ref{3.0}) holds true.
This therefore completes the proof of Lemma \ref{lem3.1}.\qed

Applying Lemma \ref{lem3.1}, we next establish the following   estimate of the unique normalized principal eigenfunction $u_\alpha$ for (\ref{1.1}) as $\alpha\rightarrow+\infty$.

\begin{lem}\label{lem3.2}
Suppose $m(x)$ satisfies (\ref{1.4}), and assume $0\leq V (x)\in C^{\gamma}(\bar{\Omega})~(0<\gamma<1)$.  Then for any fixed $\varepsilon>0$, the unique normalized principal eigenfunction $u_\alpha$ of (\ref{1.1}) satisfies
\begin{equation}\label{3.10}
w_\alpha(x):=\alpha^{-\frac{N}{4}}u_\alpha(\alpha^{-\frac{1}{2}}x)=Q(x)+o(\alpha^{-1})\ \ \text{in}\ \ \R^N\ \, \text{as}\ \ \alpha\rightarrow+\infty,
\end{equation}
where $u_\alpha(x)\equiv0$ in $\mathbb{R}^{N}\backslash\Omega$, and $Q(x)>0$ is as in (\ref{0.01}).
\end{lem}

\noindent\textbf{Proof.}
We shall complete the proof of Lemma \ref{lem3.2} by two steps.

Step 1. We claim that the function $w_\alpha(x)$ defined by (\ref{2.6}) satisfies
\begin{equation}\label{3.101}
w_\alpha(x)=\Big(\int_{\R^N}w_\alpha Qdx\Big)Q(x)+o(\alpha^{-1})\ \ \text{in}\ \ \R^N\ \ \text{as}\ \
\alpha\rightarrow+\infty,
\end{equation}
where $Q(x)>0$ is as in (\ref{0.01}).

Indeed, define
\begin{equation}\label{3.11}
f_{\alpha}(x):=w_\alpha(x)-\Big(\int_{\R^N}w_\alpha Qdx\Big)Q(x) \ \ \mbox{in}\ \, \R^N.
\end{equation}
It follows from (\ref{3.11}) that
\begin{equation}\label{3.12}
\int_{\R^N}f_{\alpha}(x)Q(x)dx\equiv0\ \ \text{for any}\ \, \alpha>0,
\end{equation}
where $\int_{\R^N}Q^2(x)dx=1$ is used.
We also define the operator $\mathcal{N}_{\alpha}$ by
\begin{equation}\label{3.13}
\mathcal{N}_{\alpha}:=-\varepsilon\Delta+\frac{4}{\varepsilon}
\sum_{i=1}^{N}a_{i}^{2}x_{i}^{2}+\alpha^{-1}V(\alpha^{-\frac{1}{2}}x)-\alpha^{-1}\lambda_{\alpha}\ \ \text{in}\ \, \Omega_{\alpha},
\end{equation}
where $\varepsilon>0$ is fixed,  $a_{i}\in\R\setminus\{0\}$  is as in (\ref{1.4}) for $i=1,\cdots,N$, and the domain $\Omega_{\alpha}$ is as in (\ref{2.6}). One can obtain from (\ref{A3}), (\ref{3.2}) and (\ref{2.7}) that for any $\alpha>0$,
\begin{equation}\label{3.14}
\mathcal{N}_{\alpha}w_\alpha=0\ \ \text{in}\ \ \Omega_{\alpha},
\end{equation}
and
\begin{equation}\label{3.15}
\mathcal{N}_{\alpha}Q=\Big[\alpha^{-1}V(\alpha^{-\frac{1}{2}}x)
-\Big(\alpha^{-1}\lambda_{\alpha}-2\sum_{i=1}^{N}|a_{i}|\Big)\Big]Q\ \,\ \text{in}\ \, \Omega_{\alpha},
\end{equation}
where the domain $\Omega_{\alpha}$ is as in (\ref{2.6}).

Applying Lemma \ref{LemmaA}, we obtain from (\ref{2.15m}) and (\ref{2.5}) that
\begin{equation}\label{3.16}
\int_{\R^N}w_\alpha(x) Q(x)dx\rightarrow1\ \ \text{as}\ \ \alpha\rightarrow+\infty.
\end{equation}
We then calculate from (\ref{0.01}), (\ref{m11}), (\ref{3.11}) and (\ref{3.13})--(\ref{3.16}) that
\begin{equation}\label{3.17}
\begin{split}
|\mathcal{N}_{\alpha}f_\alpha|=&\Big|\mathcal{N}_{\alpha}w_\alpha-\Big(\int_{\R^N}w_\alpha Q dx\Big)\mathcal{N}_{\alpha}Q\Big|\\
=&\Big|-\Big(\int_{\R^N}w_\alpha Q dx\Big)\Big[\alpha^{-1}V(\alpha^{-\frac{1}{2}}x)-
\Big(\alpha^{-1}\lambda_{\alpha}-2\sum_{i=1}^{N}|a_{i}|\Big)\Big]Q\Big|\\
=&\Big|-[1+o(1)]\cdot\Big[\alpha^{-1}V(0)+o(\alpha^{-1})-
\Big(\alpha^{-1}\lambda_{\alpha}-2\sum_{i=1}^{N}|a_{i}|\Big)\Big]Q\Big|\\
\leq&C \delta_\alpha \alpha^{-1}e^{-\frac{|x|}{4}}\ \ \text{in}\ \ B_{r_1\alpha^{\frac{1}{4}}}(0)\ \ \text{as}\ \ \alpha\rightarrow+\infty,
\end{split}
\end{equation}
where $r_1>0$ is small enough, and  $\delta_\alpha>0$ satisfies $\delta_\alpha=o(1)$ as $\alpha\rightarrow+\infty$.

Following (\ref{3.17}), we now claim that there exists a constant $C>0$, independent of $\alpha>0$, such that
\begin{equation}\label{3.18}
|f_\alpha(x)|\leq C\delta_{\alpha}\alpha^{-1}\ \ \text{uniformly in}\ \ \R^{N}\ \ \text{as}\ \ \alpha\rightarrow+\infty,
\end{equation}
where $\delta_\alpha>0$ satisfies $\delta_\alpha=o(1)$ as $\alpha\rightarrow+\infty$.
On the contrary, assume that (\ref{3.18}) is false, $i.e.$, up to a subsequence if necessary, there exists a constant $C>0$, independent of $\alpha>0$, such that
\begin{equation}\label{3.19}
\alpha\|f_\alpha\|_{L^{\infty}(\R^N)}\geq C\ \ \text{uniformly in}\ \ \R^{N}\ \ \text{as}\ \ \alpha\rightarrow+\infty.
\end{equation}
We next consider $\bar{f}_\alpha(x):=\frac{f_\alpha(x)}{\|f_\alpha\|_{L^{\infty}(\R^N)}}$, so that
\begin{equation}\label{3.20}
\|\bar{f}_{\alpha}\|_{L^{\infty}(\R^{N})}\equiv1\ \ \text{and}\ \ \int_{\R^{N}}\bar{f}_\alpha(x)Q(x)dx\equiv0\ \ \text{for any}\, \ \alpha>0
\end{equation}
in view of (\ref{3.12}). It also follows from (\ref{3.17}) and (\ref{3.19}) that there exists a constant $C>0$, independent of $\alpha>0$, such that
\begin{equation}\label{3.20m}
|\mathcal{N}_{\alpha}\bar{f}_\alpha|=\frac{|\mathcal{N}_{\alpha}f_\alpha|}{\|f_\alpha\|_{L^{\infty}(\R^{N})}}
\leq C\delta_{\alpha}e^{-\frac{|x|}{4}}\ \,\ \text{in}\ \ B_{r_1\alpha^{\frac{1}{4}}}(0)\ \ \text{as}\ \
\alpha\rightarrow+\infty,
\end{equation}
where the small constant $r_1>0$ is as in (\ref{3.17}), and  $\delta_\alpha>0$ satisfies $\delta_\alpha=o(1)$ as $\alpha\rightarrow+\infty$.

Moreover, we derive from (\ref{0.01}), (\ref{2.15m}), (\ref{3.7}), (\ref{3.11}) and (\ref{3.19}) that for above small constant $r_1>0$,
\begin{equation}\label{3.21}
\begin{split}
|\bar{f}_\alpha(x)|=&\frac{|w_{\alpha}(x)-\big(\int_{\R^N}w_{\alpha}Qdx\big)Q(x)|}
{\|f_\alpha\|_{L^{\infty}(\R^{N})}}\\[1mm]
\leq& C\alpha e^{-|x|}
\leq C|\alpha^{\frac{1}{4}}x^{-1}|^{4}e^{-\frac{|x|}{2}}\\[1mm]
\leq& C(r_1)e^{-\frac{|x|}{4}}
\ \ \text{in}\ \ \R^{N}\backslash B_{r_1\alpha^{\frac{1}{4}}}(0)\ \ \text{as}\ \
\alpha\rightarrow+\infty,
\end{split}
\end{equation}
and
\begin{equation}\label{3.22}
\begin{split}
|\nabla\bar{f}_\alpha(x)|=&\frac{|\nabla w_{\alpha}(x)-\big(\int_{\R^N}w_{\alpha}Qdx\big)\nabla Q(x)|}
{\|f_\alpha\|_{L^{\infty}(\R^{N})}}\\[1mm]
\leq& C\alpha e^{-\frac{|x|}{2}}
\leq C|\alpha^{\frac{1}{4}}x^{-1}|^{4}e^{-\frac{|x|}{3}}\\[1mm]
\leq& C(r_1)e^{-\frac{|x|}{5}}
\ \ \text{in}\ \ \R^{N}\backslash B_{r_1\alpha^{\frac{1}{4}}}(0)\ \ \text{as}\ \
\alpha\rightarrow+\infty.
\end{split}
\end{equation}
Additionally, by the definition of $\mathcal{N}_{\alpha}$ in (\ref{3.13}), we also deduce from (\ref{2.0}) and (\ref{3.20m})--(\ref{3.22}) that there exist constants $c_1>0$ and $c_2>0$, independent of $\alpha>0$, such that
\begin{equation}\label{3.23}
\begin{split}
\|\bar{f}_{\alpha}\|_{H^{1}(B_{r_1\alpha^{\frac{1}{4}}}(0))}^{2}
=&\int_{B_{r_1\alpha^{\frac{1}{4}}}(0)}|\nabla \bar{f}_{\alpha}|^{2}dx+
\int_{B_{r_1\alpha^{\frac{1}{4}}}(0)}|\bar{f}_{\alpha}|^{2}dx\\
\leq&\max\Big\{\frac{1}{\varepsilon}, 1\Big\}\int_{B_{r_1\alpha^{\frac{1}{4}}}(0)}\Big(\mathcal{N}_{\alpha}\bar{f}_{\alpha}\Big)\bar{f}_{\alpha}dx
+c_1\leq c_2\ \ \text{as}\ \ \alpha\rightarrow+\infty,
\end{split}
\end{equation}
where  the small constant $r_1>0$ is as in (\ref{3.17}).

It then follows from (\ref{3.21})--(\ref{3.23}) that the sequence $\{\bar{f}_{\alpha}\}$ is bounded uniformly in $H^{1}(\R^{N})$ as $\alpha\rightarrow+\infty$. Therefore, there exists a function $\bar{f}_{0}(x)\in H^{1}(\R^{N})$ such that, up to a subsequence if necessary, $\bar{f}_{\alp}(x)\rightarrow\bar{f}_{0}(x)$ weakly in $H^{1}(\R^{N})$ and strongly in $L_{loc}^{2}(\R^{N})$ as $\alp\rightarrow+\infty$. We thus obtain from (\ref{0.01}), (\ref{2.0}), (\ref{3.1}), (\ref{3.13}), (\ref{3.20}) and (\ref{3.20m}) that
$$\int_{\R^{N}}\bar{f}_{0}Qdx=0\ \ \text{and}\ \ \mathcal{N}\bar{f}_{0}=0\ \ \text{in}\ \ \R^{N}.$$
Since $\int_{\R^{N}}Q^2dx=1$, this further implies from  (\ref{3.4}) that
\begin{equation}\label{3.24}
\bar{f}_{0}(x)\equiv0\ \ \text{in}\ \ \R^N.
\end{equation}
On the other hand, applying the comparison principle, we obtain from (\ref{3.20m}) that there exist constants $C>0$ and $R_2>0$, independent of $\alpha>0$, such that
\begin{equation}\label{3.25}
|\bar{f}_\alpha(x)|\leq Ce^{-\frac{|x|}{6}}\ \ \text{in}\ \ B_{r_1\alpha^{\frac{1}{4}}}(0)\backslash B_{R_2}(0)\ \ \text{as}\ \ \alpha\rightarrow+\infty,
\end{equation}
where the small constant $r_1>0$ is as above.
Let $y_\alpha$ be a maximum point of $|\bar{f}_{\alpha}(x)|$ in $\R^N$, $i.e.$, $|\bar{f}_{\alpha}(y_\alpha)|=\max_{x\in\R^N}
\frac{|f_{\alpha}(x)|}{\|f_{\alpha}\|_{L^{\infty}(\R^N)}}\equiv1$.
It then follows from (\ref{3.21}) and (\ref{3.25}) that there exists a constant $C>0$, independent of $\alp>0$, such that
$$|y_\alp|\leq C\ \ \text{uniformly in}\ \ \alp>0.$$
We hence obtain that up to a subsequence if necessary, $1\equiv|\bar{f}_\alp(y_\alp)|\rightarrow|\bar{f}_{0}(\bar{y}_{0})|$ as $\alp\rightarrow+\infty$ for some $\bar{y}_0\in\R^{N}$, which however contradicts with (\ref{3.24}). Therefore, the claim (\ref{3.18}) holds true, and the claim (\ref{3.101}) is thus true in view of (\ref{3.11}).

Step 2. We claim that (\ref{3.10}) holds true.

Indeed, following (\ref{0.01}), (\ref{2.5}) and (\ref{3.101}), we can denote
\begin{equation}\label{3.26}
w_\alp(x):=Q(x)+g_\alp(x)+o(\alp^{-1})\ \ \text{in}\ \ \R^N\ \ \text{as}\ \ \alp\rightarrow+\infty,
\end{equation}
where $g_\alp(x)$ satisfies
\begin{equation}\label{3.27}
\begin{split}
\int_{\R^N}g_\alp(x) Q(x)dx\rightarrow0\ \ \text{as}\ \ \alp\rightarrow+\infty,\qquad\\
g_\alp(x)\rightarrow0\ \ \text{strongly in}\ \ L^{\infty}(\R^N)\ \ \text{as}\ \ \alp\rightarrow+\infty.
\end{split}
\end{equation}
Substituting (\ref{3.26}) into (\ref{3.101}), we obtain that
\begin{equation}\label{3.27m}
g_\alp(x)=\Big(\int_{\R^N}g_\alp Qdx\Big)Q(x)+o(\alp^{-1})\ \,\ \text{in}\ \ \R^N\ \ \text{as}\ \ \alp\rightarrow+\infty,
\end{equation}
due to the fact that  $\int_{\R^N}Q^2(x)dx=1$.

In order to complete the proof of Lemma \ref{lem3.2}, the rest is to estimate the integral $\int_{\R^N}g_\alp Qdx$. Towards this end, multiplying (\ref{3.101}) by $w_\alp(x)$ and integrating over $\R^N$, one can deduce from (\ref{0.01}), (\ref{M2.7}), (\ref{2.15m}) and (\ref{3.26}) that
\begin{equation*}
\begin{split}
1=&\int_{\R^N}w_\alp^2(x)dx=\Big(\int_{\R^N}w_\alp Qdx\Big)^2+o(\alp^{-1})\\
=&\Big[1+\int_{\R^N}g_\alp Qdx+o(\alp^{-1})\Big]^2+o(\alp^{-1})\ \ \text{as}\ \ \alp\rightarrow+\infty,
\end{split}
\end{equation*}
which implies that
\begin{equation}\label{3.28}
\Big(\int_{\R^N}g_\alp Qdx\Big)^2+2\int_{\R^N}g_\alp Qdx+o(\alp^{-1})=0\ \ \text{as}\ \ \alp\rightarrow+\infty.
\end{equation}
It then follows from (\ref{3.27}) and (\ref{3.28}) that
\begin{equation}\label{3.29}
\int_{\R^N}g_\alp(x) Q(x)dx=o(\alp^{-1})\ \ \text{as}\ \ \alp\rightarrow+\infty.
\end{equation}
We thus conclude from (\ref{3.26}), (\ref{3.27m}) and (\ref{3.29}) that (\ref{3.10}) holds true, and the proof of Lemma \ref{lem3.2} is therefore complete.\qed

\subsection{Proof of Theorem \ref{thm1.3}}
Applying Lemmas \ref{lem3.1} and \ref{lem3.2}, in this subsection we shall finish the proof of Theorem \ref{thm1.3}. The argument of this subsection can be used to obtain the higher-order asymptotic expansions of the unique principal eigenpair $(\lambda_{\alp}, u_\alpha)$ for (\ref{1.1}) as $\alp\rightarrow+\infty$, provided that $V(x)$ is smooth enough near the origin.

\vskip 0.05truein

\noindent\textbf{Proof of Theorem \ref{thm1.3}.}
We first prove that the unique principal eigenvalue $\lambda_{\alp}$ of (\ref{1.1}) satisfies
\begin{equation}\label{1.8m}
\lambda_{\alpha}=2\alpha\sum_{i=1}^{N}|a_{i}|+V(0)+\alpha^{-\frac{1}{2}}\int_{\R^N}\big[x\cdot \nabla V(0)\big]Q^2dx
+o(\alpha^{-\frac{1}{2}})\ \ \text{as}\ \ \alpha\rightarrow+\infty,
\end{equation}
where $a_{i}\in\R\setminus\{0\}$ is as in (\ref{1.4}) for $i=1,\cdots,N$, and $Q(x)>0$ is as in (\ref{0.01}).
Following (\ref{3.0}), denote
\begin{equation}\label{3.30}
\alp^{-1}\lambda_\alp:=2\sum_{i=1}^{N}|a_{i}|+\alp^{-1}V(0)+\eta_\alp+o(\eta_\alp)
\ \ \text{as}\ \ \alp\rightarrow+\infty,
\end{equation}
where $\eta_\alp$ satisfies $\lim_{\alp\rightarrow+\infty}\alp\eta_\alp=0$. Under the assumptions of Theorem \ref{thm1.3}, one can calculate from (\ref{3.6}), (\ref{3.8}), (\ref{3.10}) and (\ref{3.30}) that
\begin{equation}\label{3.30m}
\begin{split}
o(e^{-r_0\alp^{\frac{1}{4}}})
=&\int_{B_{r_0\alp^{\frac{1}{4}}}(0)}(\mathcal{N}z_\alp)Qdx\\
=&\int_{B_{r_0\alp^{\frac{1}{4}}}(0)}\Big[-\alpha^{-1}V(\alpha^{-\frac{1}{2}}x)
+\Big(\alpha^{-1}\lambda_{\alpha}-2\sum_{i=1}^{N}|a_{i}|\Big)\Big]w_\alp Qdx\\
=&\int_{B_{r_0\alp^{\frac{1}{4}}}(0)}\Big[-\alp^{-\frac{3}{2}}\big(x\cdot\nabla V(0)\big)+o(\alp^{-\frac{3}{2}}|x|)+\eta_\alp+o(\eta_\alp)\Big]\\
&\qquad\qquad\cdot\big[Q(x)+o(\alp^{-1})\big]Qdx\\
=&-\alp^{-\frac{3}{2}}\int_{\R^N}\big[x\cdot\nabla V(0)\big]Q^2dx+o(\alp^{-\frac{3}{2}})+\eta_\alp+o(\eta_\alp)\ \ \text{as}\ \ \alpha\rightarrow+\infty,
\end{split}
\end{equation}
where the operator $\mathcal{N}$ is as in (\ref{3.1}). This implies that
\begin{equation}\label{3.31}
\eta_\alp=\alp^{-\frac{3}{2}}\int_{\R^N}\big[x\cdot\nabla V(0)\big]Q^2dx+o(\alp^{-\frac{3}{2}})\ \ \text{as}\ \ \alpha\rightarrow+\infty.
\end{equation}
It then follows from (\ref{3.30}) and (\ref{3.31}) that (\ref{1.8m}) holds true.

We next claim that the unique normalized principal eigenfunction $u_\alp$ of (\ref{1.1}) satisfies
\begin{equation}\label{1.9m}
w_{\alpha}(x):=\alpha^{-\frac{N}{4}}u_{\alpha}
(\alpha^{-\frac{1}{2}}x)
=Q(x)+\alpha^{-\frac{3}{2}}\varphi_1(x)
+o(\alpha^{-\frac{3}{2}}) \ \ \text{in}\ \ \mathbb{R}^{N}
\ \ \text{as}\ \ \alpha\rightarrow+\infty,
\end{equation}
where $Q(x)>0$ is as in (\ref{0.01}), and $\varphi_1(x)$ is uniquely given by
\begin{equation}\label{3.32}
\left\{
\begin{aligned}
&-\varepsilon\Delta \varphi_{1}+\Big(\frac{4}{\varepsilon}\sum_{i=1}^{N}a_i^2x_i^2-2\sum_{i=1}^{N}|a_{i}| \Big)\varphi_{1}\\
&=\Big[\int_{\R^N}\big(x\cdot\nabla V(0)\big)Q^2dx-x\cdot\nabla V(0)\Big]Q\ \,\ \text{in}\ \,  \mathbb{R}^{N},\\[2mm]
&\int_{\mathbb{R}^{N}}\varphi_{1}Qdx=0.
\end{aligned}\right.
\end{equation}
Similar to \cite[Lemma 3.4]{Guo}, one can get the existence and uniqueness of $\varphi_1(x)$. Applying the comparison principle, it yields from (\ref{0.01}) and (\ref{3.32}) that there exists a constant $C>0$ such that
\begin{equation}\label{3.46}
|\varphi_1(x)|\leq Ce^{-\frac{|x|}{2}}\ \ \text{in}\ \, \R^N.
\end{equation}
Furthermore, it follows from (\ref{3.13}) and (\ref{3.32}) that $\varphi_1(x)$ satisfies
\begin{equation}\label{3.47}
\begin{split}
\mathcal{N}_\alp\varphi_1
=&\Big[\int_{\R^N}\big(x\cdot\nabla V(0)\big)Q^2dx-x\cdot\nabla V(0)\Big]Q\\
&+\Big[\alpha^{-1}V(\alpha^{-\frac{1}{2}}x)-\Big(\alpha^{-1}\lambda_\alpha
-2\sum_{i=1}^{N}|a_{i}|\Big)\Big]\varphi_1\
\ \text{in}\ \, \Omega_\alp,
\end{split}
\end{equation}
where the domain $\Omega_\alp$ is as in (\ref{2.6}).

Applying (\ref{3.10}) and (\ref{1.8m}), one can calculate from (\ref{0.01}), (\ref{3.14}), (\ref{3.15}), (\ref{3.46}) and (\ref{3.47}) that there exists
a constant $C>0$, independent of $\alp>0$, such that
\begin{align}
&\Big|\mathcal{N}_{\alp}\Big[w_\alp-\Big(\int_{\R^N}w_\alp Qdx\Big)Q-\alp^{-\frac{3}{2}}\varphi_1\Big]\Big|\nonumber\\
=&\Big|-[1+o(\alp^{-1})]\mathcal{N}_{\alp}Q-\alp^{-\frac{3}{2}}\mathcal{N}_{\alp}\varphi_1\Big|\nonumber\\
=&\Big|-[1+o(\alp^{-1})]\cdot\Big[\alpha^{-1}V(\alpha^{-\frac{1}{2}}x)-\Big(
\alp^{-\frac{3}{2}}\int_{\R^N}\big(x\cdot\nabla V(0)\big)Q^2dx\label{3.48}\\
&+\alp^{-1}V(0)+o(\alp^{-\frac{3}{2}})
\Big)\Big]Q-\alp^{-\frac{3}{2}}\Big[\alpha^{-1}V(\alpha^{-\frac{1}{2}}x)
-\Big(\alpha^{-1}\lambda_\alpha-2\sum_{i=1}^{N}|a_{i}|\Big)\Big]\varphi_1\nonumber\\
&-\alp^{-\frac{3}{2}}\Big[\int_{\R^N}\big(x\cdot\nabla V(0)\big)Q^2dx-x\cdot\nabla V(0)\Big]Q\Big|\nonumber\\
\leq& C\delta_\alp \alp^{-\frac{3}{2}}e^{-\frac{|x|}{4}}\ \ \text{in}\ \ B_{r_{0}\alpha^{\frac{1}{4}}}(0)\ \ \text{as}\ \ \alpha\rightarrow+\infty,\nonumber
\end{align}
where $r_0>0$ is given by Theorem \ref{thm1.3}, and $\delta_\alp>0$ satisfies $\delta_\alpha=o(1)$ as $\alpha\rightarrow+\infty$. Similar to Step 1 in the proof of Lemma \ref{lem3.2}, it yields from (\ref{3.48}) that
\begin{equation}\label{3.49}
w_\alp(x)=\Big(\int_{\R^N}w_\alp Qdx\Big)Q(x)+\alp^{-\frac{3}{2}}\varphi_1(x)+o(\alp^{-\frac{3}{2}})
\ \ \text{in}\ \ \R^N\ \,\ \text{as}\ \ \alp\rightarrow+\infty.
\end{equation}
Moreover, similar to Step 2 in the proof of Lemma \ref{lem3.2}, we further obtain from (\ref{3.49}) that the claim (\ref{1.9m}) holds true.

We now prove that (\ref{1.8}) holds true. Following (\ref{1.8m}), denote
\begin{equation}\label{3.50}
\begin{split}
\alp^{-1}\lambda_\alp:=&2\sum_{i=1}^{N}|a_{i}|+\alp^{-1}V(0)+\alp^{-\frac{3}{2}}\int_{\R^N}\big[x\cdot\nabla V(0)\big]Q^2dx\\
&+\hat{\eta}_\alp+o(\hat{\eta}_\alp)
\ \ \text{as}\ \ \alp\rightarrow+\infty,
\end{split}
\end{equation}
where $\hat{\eta}_{\alpha}$ satisfies $\lim_{\alpha\rightarrow+\infty}\alpha^{\frac{3}{2}}\hat{\eta}_{\alpha}=0$. It then follows from (\ref{3.30m}), (\ref{1.9m}) and (\ref{3.50}) that
\begin{align}
o(e^{-r_0\alp^{\frac{1}{4}}})
=&\int_{B_{r_0\alp^{\frac{1}{4}}}(0)}(\mathcal{N}z_\alp)Qdx\nonumber\\
=&\int_{B_{r_0\alp^{\frac{1}{4}}}(0)}\Big[-\alpha^{-1}V(\alpha^{-\frac{1}{2}}x)
+\Big(\alpha^{-1}\lambda_{\alpha}-2\sum_{i=1}^{N}|a_{i}|\Big)\Big]w_\alp Qdx\nonumber\\
=&\int_{B_{r_0\alp^{\frac{1}{4}}}(0)}\Big[-\alp^{-\frac{3}{2}}\big(x\cdot\nabla V(0)\big)-\alp^{-2}\sum_{|\tau|=2}\frac{D^{\tau}V(0)}{\tau!}x^{\tau}
-o(\alp^{-2}|x|^{2})\label{3.51}\\
&\qquad\qquad+\alp^{-\frac{3}{2}}\int_{\R^N}\big(x\cdot\nabla V(0)\big)Q^2dx+\hat{\eta}_\alp+o(\hat{\eta}_\alp)\Big]\nonumber\\
&\qquad\qquad\cdot\big[Q(x)+\alpha^{-\frac{3}{2}}\varphi_1(x)
+o(\alpha^{-\frac{3}{2}})\big]Qdx\nonumber\\
=&-\alp^{-2}\sum_{|\tau|=2}\frac{D^{\tau}V(0)}{\tau!}\int_{\R^N}x^{\tau}Q^2dx+o(\alp^{-2})
+\hat{\eta}_\alp+o(\hat{\eta}_\alp)
\ \ \text{as}\ \ \alpha\rightarrow+\infty,\nonumber
\end{align}
which implies that
\begin{equation}\label{3.52}
\hat{\eta}_\alp=\alp^{-2}\sum_{|\tau|=2}\frac{D^{\tau}V(0)}{\tau!}\int_{\R^N}x^{\tau}Q^2dx+o(\alp^{-2})
\ \ \text{as}\ \ \alpha\rightarrow+\infty.
\end{equation}
It thus yields from (\ref{3.50}) and (\ref{3.52}) that (\ref{1.8}) holds true.

We finally prove that (\ref{1.9}) holds true. Note from (\ref{1.7}), (\ref{1.7m}) and (\ref{3.13}) that $\varphi_2(x)$ satisfies
\begin{equation}\label{3.53}
\begin{split}
\mathcal{N}_\alp\varphi_2=&\sum_{|\tau|=2}\frac{D^{\tau}V(0)}{\tau!}
\Big(\int_{\R^N}x^{\tau}Q^2dx-x^{\tau}\Big)Q
\\
&+\Big[\alpha^{-1}V(\alpha^{-\frac{1}{2}}x)
-\Big(\alpha^{-1}\lambda_\alpha-2\sum_{i=1}^{N}|a_{i}|\Big)\Big]\varphi_2\ \,\ \text{in}\ \ \Omega_\alp,
\end{split}
\end{equation}
where the domain $\Omega_\alp$ is as in (\ref{2.6}). Moreover, applying the comparison principle to (\ref{1.7}) and (\ref{1.7m}), we derive from (\ref{0.01}) that there exists a constant $C>0$ such that
\begin{equation}\label{3.54}
|\varphi_2(x)|\leq Ce^{-\frac{|x|}{2}}\ \,\ \text{in}\ \ \R^N.
\end{equation}
Applying (\ref{1.8}) and (\ref{1.9m}), we then calculate from (\ref{0.01}), (\ref{3.14}), (\ref{3.15}), (\ref{3.46}), (\ref{3.47}), (\ref{3.53}) and (\ref{3.54}) that there exists a constant $C>0$, independent of $\alp>0$, such that
\begin{align}
&\Big|\mathcal{N}_{\alp}\Big[w_\alp-\Big(\int_{\R^N}w_\alp Qdx\Big)Q-\alp^{-\frac{3}{2}}\varphi_1-\alp^{-2}\varphi_2\Big]\Big|\nonumber\\
=&\Big|\big[1+o(\alp^{-\frac{3}{2}})\big]\mathcal{N}_{\alp}Q+\alp^{-\frac{3}{2}}\mathcal{N}_{\alp}\varphi_1
+\alp^{-2}\mathcal{N}_{\alp}\varphi_2\Big|\nonumber\\
=&\Big|\big[1+o(\alp^{-\frac{3}{2}})\big]\cdot\Big[\alpha^{-1}V(\alpha^{-\frac{1}{2}}x)
-\Big(\alpha^{-1}\lambda_{\alpha}-2\sum_{i=1}^{N}|a_{i}|\Big)\Big]Q+\alp^{-\frac{3}{2}}\Big\{\Big[\alpha^{-1}V(\alpha^{-\frac{1}{2}}x)
\nonumber\\
&-\Big(\alpha^{-1}\lambda_{\alpha}-2\sum_{i=1}^{N}|a_{i}|\Big)\Big]\varphi_1+\Big[\int_{\R^N}\big(x\cdot\nabla V(0)\big)Q^2dx-x\cdot\nabla V(0)\Big]Q\Big\}\nonumber\\
&+\alp^{-2}\Big\{\Big[\alpha^{-1}V(\alpha^{-\frac{1}{2}}x)
-\Big(\alpha^{-1}\lambda_{\alpha}-2\sum_{i=1}^{N}|a_{i}|\Big)\Big]\varphi_2\nonumber\\
&+\sum_{|\tau|=2}\frac{D^{\tau}V(0)}{\tau!}
\Big(\int_{\R^N}x^{\tau}Q^2dx-x^{\tau}\Big)Q\Big\}\Big|\label{3.55}\\
\leq&\Big|\big[1+o(\alp^{-\frac{3}{2}})\big]\cdot\Big[\alp^{-\frac{3}{2}}\big(x\cdot\nabla V(0)\big)+\alp^{-2}\sum_{|\tau|=2}\frac{D^{\tau}V(0)}{\tau!}x^{\tau}
+o(\alp^{-2}|x|^{2})\nonumber\\
&-\alp^{-\frac{3}{2}}\int_{\R^N}\big(x\cdot\nabla V(0)\big)Q^2dx-\alp^{-2}\sum_{|\tau|=2}\frac{D^{\tau}V(0)}{\tau!}\int_{\R^N}x^{\tau}Q^2dx
-o(\alp^{-2})\Big]Q\nonumber\\
&+\alp^{-\frac{3}{2}}\Big[\int_{\R^N}\big(x\cdot\nabla V(0)\big)Q^2dx-x\cdot\nabla V(0)\Big]Q\nonumber\\
&+\alp^{-2}\sum_{|\tau|=2}\frac{D^{\tau}V(0)}{\tau!}
\Big(\int_{\R^N}x^{\tau}Q^2dx-x^{\tau}\Big)Q\Big|+o(\alp^{-2})\big(|\varphi_1|+|\varphi_2|\big)
\nonumber\\
\leq& C\delta_\alp \alp^{-2}e^{-\frac{|x|}{4}}\ \,\ \text{in}\ \ B_{r_{0}\alpha^{\frac{1}{4}}}(0)\ \ \text{as}\ \ \alpha\rightarrow+\infty,\nonumber
\end{align}
where $r_0>0$ is as in Theorem \ref{thm1.3}, and $\delta_\alp>0$ satisfies $\delta_\alpha=o(1)$ as $\alpha\rightarrow+\infty$.
Similar to Step 1 in the proof of  Lemma \ref{lem3.2}, we then obtain from (\ref{3.55}) that
\begin{equation}\label{3.56}
\begin{split}
w_\alp(x)=&\Big(\int_{\R^N}w_\alp Qdx\Big)Q(x)+\alp^{-\frac{3}{2}}\varphi_1(x)\\
&+\alp^{-2}\varphi_2(x)+o(\alp^{-2})
\ \ \text{in}\ \ \R^N\ \ \text{as}\ \ \alp\rightarrow+\infty.
\end{split}
\end{equation}
Finally, similar to Step 2 in the proof of Lemma \ref{lem3.2}, one can deduce from (\ref{3.56}) that (\ref{1.9}) holds true.
The proof of Theorem \ref{thm1.3} is therefore complete.
\qed

\section{Proof of Theorem \ref{cor1.3}}
Suppose the regularity of the potential $V(x)$ near the origin is lower than that of Theorem \ref{thm1.3}, the purpose of this section is to prove Theorem \ref{cor1.3} on the refined limiting profiles of $(\lambda_{\alpha}, u_{\alpha})$ for  (\ref{1.1}) as $\alp\to +\infty$. Towards this aim, we first analyze the following second-order asymptotic expansions as $\alpha\rightarrow+\infty$.

\begin{lem}\label{lem4.1}
Suppose   $m(x)$ satisfies (\ref{1.4}), and assume  $0\leq V(x)\in C^{\gamma}(\bar{\Omega})~(0<\gamma<1)$ satisfies $V(x)=h_0(x)+o(|x|^{k_0})$ as $|x|\rightarrow0$, where $h_0(x)$ satisfies (\ref{1.12}) for some $k_0>0$. Then for any fixed $\varepsilon>0$, the unique principal eigenpair $(\lambda_{\alpha}, u_{\alpha})$ of (\ref{1.1}) satisfies
\begin{equation}\label{4.1}
\lambda_{\alpha}=2\alpha\sum_{i=1}^{N}|a_{i}|+\alpha^{-\frac{k_0}{2}}\int_{\R^{N}}h_0(x)Q^2dx
+o(\alpha^{-\frac{k_0}{2}})\ \ \text{as}\ \ \alpha\rightarrow+\infty,
\end{equation}
and
\begin{equation}\label{4.2}
\begin{split}
w_{\alpha}(x):=&\alpha^{-\frac{N}{4}}u_{\alpha}
(\alpha^{-\frac{1}{2}}x)\\
=&Q(x)+\alpha^{-\frac{k_0+2}{2}}\varphi_3(x)
+o(\alpha^{-\frac{k_0+2}{2}}) \ \ \text{in}\ \ \mathbb{R}^{N}
\ \ \text{as}\ \ \alpha\rightarrow+\infty,
\end{split}
\end{equation}
where  $a_{i}\in \R\setminus\{0\}$ is as in (\ref{1.4}) for $i=1,\cdots,N$, $Q(x)>0$ is as in (\ref{0.01}), $u_\alpha(x)\equiv0$ in $\mathbb{R}^{N}\backslash\Omega$, and $\varphi_3(x)$ is given uniquely by (\ref{1.10}) and (\ref{1.11}).
\end{lem}

\noindent\textbf{Proof.}
Recall from (\ref{3.6}) that
\begin{equation}\label{4.3}
\mathcal{N}z_{\alpha}=\Big[-\alpha^{-1}V(\alpha^{-\frac{1}{2}}x)
+\Big(\alpha^{-1}\lambda_\alpha-2\sum_{i=1}^{N}|a_{i}|\Big)\Big]w_\alpha\ \ \text{in}\ \ \Omega_\alpha,
\end{equation}
where $\Omega_\alpha$, $\mathcal{N}$ and $z_\alpha$ are as in (\ref{2.6}), (\ref{3.1}) and (\ref{3.5}), respectively. Similar to (\ref{3.8}), we calculate from (\ref{0.01}), (\ref{2.15m}), (\ref{2.5}) and (\ref{4.3}) that
\begin{equation}\label{4.3m}
\begin{split}
o(e^{-\hat{r}_0\alpha^{\frac{1}{4}}})=&\int_{B_{\hat{r}_{0}\alpha^{\frac{1}{4}}}(0)}
\big(\mathcal{N}z_{\alpha}\big)
Qdx\\
=&\int_{B_{\hat{r}_{0}\alpha^{\frac{1}{4}}}(0)}\Big[-\alpha^{-1}V(\alpha^{-\frac{1}{2}}x)
+\Big(\alpha^{-1}\lambda_{\alpha}-2\sum_{i=1}^{N}|a_{i}|\Big)\Big]w_\alpha Qdx\\
=&-\alpha^{-\frac{k_0+2}{2}}\int_{\R^N}h_0(x)Q^2dx-o(\alpha^{-\frac{k_0+2}{2}})\\
&+\Big(\alpha^{-1}\lambda_{\alpha}-2\sum_{i=1}^{N}|a_{i}|\Big)[1+o(1)]\ \ \text{as}\ \ \alpha\rightarrow+\infty,
\end{split}
\end{equation}
where $\hat{r}_0>0$ is as in Theorem \ref{cor1.3}. This implies that
\begin{equation}\label{4.4}
\begin{split}
\alpha^{-1}\lambda_{\alpha}-2\sum_{i=1}^{N}|a_{i}|=\alpha^{-\frac{k_0+2}{2}}\int_{\R^N}h_0(x)Q^2dx
+o(\alpha^{-\frac{k_0+2}{2}})\ \ \text{as}\ \ \alpha\rightarrow+\infty.
\end{split}
\end{equation}
Thus, the unique principal eigenvalue $\lambda_{\alpha}$ of (\ref{1.1}) satisfies (\ref{4.1}) in view of (\ref{4.4}).

We next prove the second-order asymptotic expansion of the unique normalized principal eigenfunction $u_\alpha$ for (\ref{1.1}) as $\alpha\rightarrow+\infty$. It follows from (\ref{1.10}), (\ref{1.11}) and (\ref{3.13}) that $\varphi_3(x)$ satisfies
\begin{equation}\label{4.5}
\begin{split}
\mathcal{N}_{\alpha}\varphi_{3}=&\Big[\alpha^{-1}V(\alpha^{-\frac{1}{2}}x)-
\Big(\alpha^{-1}\lambda_{\alpha}-2\sum_{i=1}^{N}|a_{i}|\Big)\Big]
\varphi_3\\
&+\Big[\int_{\R^N}h_0(x)Q^2dx-h_0(x)\Big]Q
\ \ \text{in}\ \ \Omega_{\alpha},
\end{split}
\end{equation}
where the domain $\Omega_{\alpha}$ is as in (\ref{2.6}). Moreover, applying the comparison principle to (\ref{1.10}) and (\ref{1.11}), we derive from (\ref{0.01}) that there exists a constant $C>0$ such that
\begin{equation}\label{4.6}
|\varphi_{3}(x)|\leq Ce^{-\frac{|x|}{2}}\ \,\ \text{and}\,\ \ |\varphi_{4}(x)|\leq Ce^{-\frac{|x|}{3}}\ \ \text{in}\ \ \R^N.
\end{equation}
Under the assumptions of Lemma \ref{lem4.1}, we calculate from (\ref{0.01}), (\ref{2.5}), (\ref{3.14}), (\ref{3.15}) and (\ref{4.4})--(\ref{4.6}) that there exists a constant $C>0$, independent of $\alpha>0$, such that
\begin{align}
&\Big|\mathcal{N}_{\alpha}\Big[w_\alpha-\Big(\int_{\R^N}w_\alpha Qdx\Big)Q-
\alpha^{-\frac{k_0+2}{2}}\varphi_{3}\Big]\Big|\nonumber\\
=&\Big|[1+o(1)]\cdot\Big[\alpha^{-1}V(\alpha^{-\frac{1}{2}}x)-
\Big(\alpha^{-1}\lambda_{\alpha}-2\sum_{i=1}^{N}|a_{i}|\Big)\Big]Q+\alpha^{-\frac{k_0+2}{2}}
\Big\{\Big[\alpha^{-1}V(\alpha^{-\frac{1}{2}}x)\nonumber\\
&-
\Big(\alpha^{-1}\lambda_{\alpha}-2\sum_{i=1}^{N}|a_{i}|\Big)\Big]\varphi_{3}
+\Big(\int_{\R^N}h_0(x)Q^2dx-h_0(x)\Big)Q\Big\}\Big|\nonumber\\
=&\Big|[1+o(1)]\cdot\Big[\alpha^{-\frac{k_0+2}{2}}h_0(x)+o(\alpha^{-\frac{k_0+2}{2}}|x|^{k_0})-
\alpha^{-\frac{k_0+2}{2}}\int_{\R^N}h_0(x)Q^2dx\label{4.7}\\
&-o(\alpha^{-\frac{k_0+2}{2}})\Big]Q+\alpha^{-\frac{k_0+2}{2}}
\Big[\alpha^{-1}V(\alpha^{-\frac{1}{2}}x)-
\Big(\alpha^{-1}\lambda_{\alpha}-2\sum_{i=1}^{N}|a_{i}|\Big)\Big]\varphi_{3}\nonumber\\
&+\alpha^{-\frac{k_0+2}{2}}\Big(\int_{\R^N}h_0(x)Q^2dx-h_0(x)\Big)Q\Big|\nonumber\\
\leq&C\delta_{\alpha}\alpha^{-\frac{k_0+2}{2}}e^{-\frac{|x|}{4}}
\ \ \text{in}\ \ B_{\hat{r}_{0}\alpha^{\frac{1}{4}}}(0)\ \ \text{as}\ \ \alpha\rightarrow+\infty,\nonumber
\end{align}
where $\delta_{\alpha}>0$ satisfies $\delta_{\alpha}=o(1)$ as $\alpha\rightarrow+\infty$, and $\hat{r}_{0}>0$ is as in Theorem \ref{cor1.3}. Similar to Step 1 in the proof of Lemma \ref{lem3.2}, one can derive from (\ref{4.7}) that
\begin{equation}\label{4.8}
\begin{split}
w_\alpha(x)=\Big(\int_{\R^N}w_\alpha Qdx\Big)Q(x)+\alpha^{-\frac{k_0+2}{2}}\varphi_{3}(x)+o(\alpha^{-\frac{k_0+2}{2}})\ \ \text{in}\ \ \R^N\ \ \text{as}\ \ \alpha\rightarrow+\infty.
\end{split}
\end{equation}
Finally, similar to Step 2 in the proof of Lemma \ref{lem3.2}, it yields from (\ref{4.8}) that (\ref{4.2}) holds true. This completes the proof of Lemma \ref{4.1}.\qed

Applying Lemma \ref{4.1}, we are ready to finish the proof of Theorem \ref{cor1.3}.
\vskip 0.05truein

\noindent\textbf{Proof of Theorem \ref{cor1.3}.}
We first prove that the unique principal eigenvalue $\lambda_{\alpha}$ of (\ref{1.1}) satisfies (\ref{1.13}). In view of  (\ref{4.1}), we define
\begin{equation}\label{4.9}
\alpha^{-1}\lambda_{\alpha}-2\sum_{i=1}^{N}|a_{i}|
:=\alpha^{-\frac{k_0+2}{2}}\int_{\R^N}h_0(x)Q^2dx+\kappa_{\alpha}
+o(\kappa_{\alpha})\ \ \text{as}\ \ \alpha\rightarrow+\infty,
\end{equation}
where $\kappa_{\alpha}$ satisfies $\lim_{\alpha\rightarrow+\infty}\alpha^{\frac{k_0+2}{2}}\kappa_{\alpha}=0$.
Under the assumptions of Theorem \ref{cor1.3}, substituting (\ref{4.9}) into (\ref{4.3m}), we derive from (\ref{0.01}) and (\ref{4.2}) that
\begin{align}
o(e^{-\hat{r}_0\alpha^{\frac{1}{4}}})=&\int_{B_{\hat{r}_{0}\alpha^{\frac{1}{4}}}(0)}
\big(\mathcal{N}z_{\alpha}\big)Qdx\nonumber\\
=&\int_{B_{\hat{r}_{0}\alpha^{\frac{1}{4}}}(0)}\Big[-\alpha^{-1}V(\alpha^{-\frac{1}{2}}x)
+\Big(\alpha^{-1}\lambda_{\alpha}-2\sum_{i=1}^{N}|a_{i}|\Big)\Big]w_\alpha Qdx\nonumber\\
=&\int_{B_{\hat{r}_{0}\alpha^{\frac{1}{4}}}(0)}\Big[-\alpha^{-\frac{k_0+2}{2}}h_0(x)
+\alpha^{-\frac{k_0+2}{2}}\int_{\R^N}h_0(x)Q^2dx+\kappa_{\alpha}
+o(\kappa_{\alpha})\Big]\label{4.9m}\\
&\qquad\quad\quad\cdot\Big[Q+\alpha^{-\frac{k_0+2}{2}}\varphi_3
+o(\alpha^{-\frac{k_0+2}{2}})\Big]Qdx\nonumber\\
=&-\alpha^{-k_0-2}\int_{\R^N}h_0(x)\varphi_{3}Qdx+o(\alpha^{-k_0-2})+\kappa_{\alpha}
+o(\kappa_{\alpha})\ \ \text{as}\ \ \alpha\rightarrow+\infty,\nonumber
\end{align}
where $\hat{r}_0>0$ is as in Theorem \ref{cor1.3}. This implies that
\begin{equation}\label{4.10}
\kappa_{\alpha}=\alpha^{-k_0-2}\int_{\R^N}h_0(x)\varphi_{3}Qdx+o(\alpha^{-k_0-2})\ \ \text{as}\ \ \alpha\rightarrow+\infty.
\end{equation}
We thus conclude from (\ref{4.9}) and (\ref{4.10}) that (\ref{1.13}) holds true.

We next prove that the unique normalized principal eigenfunction $u_\alpha$ of (\ref{1.1}) satisfies (\ref{1.14}). It yields from (\ref{1.10}), (\ref{1.11}) and (\ref{3.13}) that $\varphi_4(x)$ satisfies
\begin{equation}\label{4.11}
\begin{split}
\mathcal{N}_{\alpha}\varphi_{4}=&\Big(\int_{\R^N}h_0(x)\varphi_3Qdx\Big)Q+\Big(\int_{\R^N}h_0(x)Q^{2}dx
-h_0(x)\Big)\varphi_3 \\
&+\Big[\alpha^{-1}V(\alpha^{-\frac{1}{2}}x)-
\Big(\alpha^{-1}\lambda_{\alpha}-2\sum_{i=1}^{N}|a_{i}|\Big)\Big]
\varphi_4
\ \ \text{in}\ \ \Omega_{\alpha},
\end{split}
\end{equation}
where the domain $\Omega_{\alpha}$ is as in (\ref{2.6}). Applying (\ref{1.13}) and (\ref{4.2}), we then calculate from (\ref{0.01}), (\ref{3.14}), (\ref{3.15}), (\ref{4.5}), (\ref{4.6}) and (\ref{4.11}) that there exists a constant $C>0$, independent of $\alp>0$, such that
\begin{align}
&\Big|\mathcal{N}_{\alpha}\Big[w_\alpha-\Big(\int_{\R^N}w_\alpha Qdx\Big)Q-
\alpha^{-\frac{k_0+2}{2}}\varphi_{3}-\alpha^{-k_0-2}\varphi_{4}\Big]\Big|\nonumber\\
=&\Big|[1+o(\alpha^{-\frac{k_0+2}{2}})]\cdot\Big[\alpha^{-1}V(\alpha^{-\frac{1}{2}}x)-
\Big(\alpha^{-1}\lambda_{\alpha}-2\sum_{i=1}^{N}|a_{i}|\Big)\Big]Q\nonumber\\
&+\alpha^{-\frac{k_0+2}{2}}
\Big\{\Big[\alpha^{-1}V(\alpha^{-\frac{1}{2}}x)-
\Big(\alpha^{-1}\lambda_{\alpha}-2\sum_{i=1}^{N}|a_{i}|\Big)\Big]\varphi_{3}
+\Big[\int_{\R^N}h_0(x)Q^2dx-h_0(x)\Big]Q\Big\}\nonumber\\
&+\alpha^{-k_0-2}
\Big\{\Big(\int_{\R^N}h_0(x)\varphi_3Qdx\Big)Q+\Big(\int_{\R^N}h_0(x)Q^{2}dx
-h_0(x)\Big)\varphi_3
\label{4.11m}\\
&+\Big[\alpha^{-1}V(\alpha^{-\frac{1}{2}}x)-
\Big(\alpha^{-1}\lambda_{\alpha}-2\sum_{i=1}^{N}|a_{i}|\Big)\Big]
\varphi_4\Big\}\Big|\nonumber\\
=&\Big|o(\alpha^{-k_0-2})\big[\big(1+h_0(x)\big)Q+\varphi_{3}\big]\nonumber\\
&+\alpha^{-k_0-2}
\Big[\alpha^{-1}V(\alpha^{-\frac{1}{2}}x)-
\Big(\alpha^{-1}\lambda_{\alpha}-2\sum_{i=1}^{N}|a_{i}|\Big)\Big]
\varphi_4\Big|\nonumber\\
\leq&C\delta_{\alpha}\alpha^{-k_0-2}e^{-\frac{|x|}{4}}
\ \ \text{in}\ \ B_{\hat{r}_{0}\alpha^{\frac{1}{4}}}(0)\ \ \text{as}\ \ \alpha\rightarrow+\infty,\nonumber
\end{align}
where $\hat{r}_{0}>0$ is as in Theorem \ref{cor1.3}, and $\delta_{\alpha}>0$ satisfies $\delta_{\alpha}=o(1)$ as $\alpha\rightarrow+\infty$.
Similar to Step 1 in the proof of Lemma \ref{lem3.2}, it then yields from (\ref{4.11m}) that
\begin{equation}\label{4.12}
\begin{split}
w_\alpha(x)=&\Big(\int_{\R^N}w_\alpha Qdx\Big)Q(x)+\alpha^{-\frac{k_0+2}{2}}\varphi_{3}(x)+\alpha^{-k_0-2}\varphi_{4}(x)\\
&+o(\alpha^{-k_0-2})\ \ \text{in}\ \ \R^N\ \ \text{as}\ \ \alpha\rightarrow+\infty.
\end{split}
\end{equation}

Following (\ref{4.2}) and (\ref{4.12}), we define
\begin{equation}\label{4.13}
w_{\alpha}(x):=Q(x)+\alpha^{-\frac{k_0+2}{2}}\varphi_{3}(x)+\beta_{\alpha}(x)+o(\alpha^{-k_0-2})
\ \ \text{in}\ \ \R^N\ \ \text{as}\ \ \alpha\rightarrow+\infty,
\end{equation}
where $\beta_{\alpha}(x)$ satisfies
\begin{equation}\label{4.14}
\begin{split}
\alpha^{\frac{k_0+2}{2}}\int_{\R^N}\beta_{\alpha}(x) Q(x)dx\rightarrow0\ \ \text{as}\ \ \alp\rightarrow+\infty,\qquad\\
\alpha^{\frac{k_0+2}{2}}\beta_{\alpha}(x)\rightarrow0\ \ \text{strongly in}\ \ L^{\infty}(\R^N)\ \ \text{as}\ \ \alp\rightarrow+\infty.
\end{split}
\end{equation}
Multiplying (\ref{4.12}) by $w_{\alpha}(x)$ and integrating over $\R^N$, we obtain from (\ref{0.01}), (\ref{M2.7}), (\ref{2.15m}), (\ref{4.6}), (\ref{4.13}) and (\ref{4.14}) that
\begin{equation}\label{4.14m}
\begin{split}
1=&\int_{\R^N}w_{\alpha}^2(x)dx\\
=&\Big(\int_{\R^N}w_{\alpha}Qdx\Big)^{2}
+\alpha^{-k_0-2}\int_{\R^N}\varphi_{3}^2(x)dx+o(\alpha^{-k_0-2})\\
=&\Big(1+\int_{\R^N}\beta_{\alpha}(x)Qdx+o(\alpha^{-k_0-2})\Big)^{2}
+\alpha^{-k_0-2}\int_{\R^N}\varphi_{3}^2(x)dx+o(\alpha^{-k_0-2})\\
=&1+2\int_{\R^N}\beta_{\alpha}(x)Qdx
+\alpha^{-k_0-2}\int_{\R^N}\varphi_{3}^2(x)dx+o(\alpha^{-k_0-2})\ \ \text{as}\ \ \alp\rightarrow+\infty,
\end{split}
\end{equation}
where we have used the facts that
$$\int_{\R^N}Q^2(x)dx=1\ \ \text{and}\ \ \int_{\R^N}\varphi_3Qdx=\int_{\R^N}\varphi_4Qdx=0.$$
It then follows from (\ref{4.14m}) that
\begin{equation}\label{4.15}
\int_{\R^N}\beta_{\alpha}(x)Qdx=-\frac{\alpha^{-k_0-2}}{2}\int_{\R^N}\varphi_{3}^2(x)dx
+o(\alpha^{-k_0-2})\ \ \text{as}\ \ \alp\rightarrow+\infty.
\end{equation}
Additionally, substituting (\ref{4.13}) into (\ref{4.12}), we obtain from (\ref{4.15}) that
\begin{equation}\label{4.16}
\begin{split}
\beta_{\alpha}(x)=&\Big(\int_{\R^N}\beta_{\alpha}(x)Qdx\Big)Q(x)
+\alpha^{-k_0-2}\varphi_{4}(x)+o(\alpha^{-k_0-2})\\
=&\alpha^{-k_0-2}\Big[\varphi_{4}(x)-\frac{1}{2}\Big(\int_{\R^N}\varphi_{3}^2(x)dx\Big)Q(x)\Big]\\
&+o(\alpha^{-k_0-2})\ \,\ \text{in}\ \ \R^N\ \ \text{as}\ \ \alp\rightarrow+\infty.
\end{split}
\end{equation}
We finally conclude from (\ref{4.13}) and (\ref{4.16}) that (\ref{1.14}) holds true. This therefore completes the proof of Theorem \ref{cor1.3}.\qed

\vskip 0.4truein
\noindent {\bf Acknowledgements}
\vskip 0.1truein
Yujin Guo is partially supported by National Key R $\&$ D Program of China (Grant 2023YFA1010001), and NSF of China (Grants 12225106 and 12371113). Yuan Lou is partially supported by NSF of China (Grants
12250710674 and 12261160366).

\vskip 0.3truein

\noindent {\bf Data availability}
\vskip 0.1truein
No data was used for the research described in the article.

\vskip 0.3truein
\noindent {\bf Declaration of competing interest}
\vskip 0.1truein
The authors declare that there is no conflict of interests regarding the publication of this paper.

\end{document}